\documentclass[10pt,a4paper,english]{article} 

\usepackage{color}
\usepackage[utf8]{inputenc}

\usepackage{amsmath, amssymb, theorem, latexsym}
\usepackage[english]{babel}

\usepackage{float}

\usepackage{stmaryrd}
\usepackage{enumerate}
\usepackage{eufrak}
\usepackage{graphicx}

\usepackage[margin=1in]{geometry}

\def\CA{{\mathcal A}}

\def\CD{{\mathcal D}}
\def\CE{{\mathcal E}}

\def\CH{{\mathcal H}}

\def\CL{{\mathcal L}}

\def\CR{{\mathcal R}}

\def\Hess{{\rm Hess}}

\def\pa{{\partial}}

\def\curl{ \mbox{ curl }}

\newcommand{\sca}[2]{\langle#1,#2\rangle}
\newcommand{\barr}{\overline}

\newcommand{\beq}{\begin{equation}}
\newcommand{\eeq}{\end{equation}}

\newcommand{\Supp}{\textrm{Supp~}}

\renewcommand{\Re}{{\rm Re\,}}
\renewcommand{\Im}{{\rm Im\,}}
\newtheorem{theorem}{Theorem}[section]
\newtheorem{lemma}[theorem]{Lemma}

\newtheorem{proposition}[theorem]{Proposition}
\newtheorem{definition}[theorem]{Definition}

\newtheorem{corollary}[theorem]{Corollary}

\makeatletter
\@addtoreset{equation}{section}

\makeatother

\title{On the semi-classical analysis of Schr\"odinger operators with purely
imaginary electric potentials in a bounded domain}

\author{Raphaël \textsc{Henry}\footnote{
D\'epartement de Math\'ematiques, Batiment 425, Universit\'e Paris Sud, 91405 Orsay Cedex, France.
email: raphael.henry@math.u-psud.fr}
\thanks{The author is supported  by the ANR NOSEVOL.}}

\date{}

\begin{document}
\maketitle

\begin{abstract}
In this paper, we describe the leftmost eigenvalue of the non-selfadjoint operator $\CA_h = -h^2\Delta+iV(x)$ 
 with Dirichlet boundary conditions on a smooth bounded domain $\Omega\subset\mathbb{R}^n\,$,
 as $h\rightarrow0\,$. $V$ is assumed to be a Morse function without critical point at the boundary of $\Omega\,$. 
 More precisely, we compare $\inf\Re\sigma(\CA_h)$ with the minimum of the spectrum's real part for some model operator. In
 the case where $V$ has no critical point, the spectrum is determined by the boundary points where $\nabla V$ is orthogonal, and
 the model operator involves a $1$-dimensional complex Airy operator in $\mathbb{R}^+\,$. If $V$ is a Morse function with critical points in
 $\Omega\,$,
 the behavior of the operator near the critical points prevails, and the model operator is a complex harmonic oscillator.\\
 This question is related to the decay of associated semigroups. In particular, it allows to recover, in a simplified setting, some stability 
 results of \cite{Alm} in superconductivity theory.
\end{abstract}

\section{Introduction}
Let $n\geq1$, $h_0>0$, and $\Omega\subset\mathbb{R}^n$ be a smooth bounded domain. We consider, for $h\in(0,h_0)$, the operator
\beq\label{defAh}
\CA_h = -h^2\Delta+iV(x)\,,~~~\CD(\CA_h) = H_0^1(\Omega ; \mathbb{C})\cap H^2(\Omega ; \mathbb{C})\,,
\eeq
where $V\in\mathcal{C}^\infty(\Omega\ ; \mathbb{R})$ is a smooth potential.\\
Under these conditions $\CA_h$ has compact resolvent, hence discrete spectrum and the purpose of this paper is to understand the behavior as $h\rightarrow0$ of the smallest real part of $\lambda(h)$, for
$\lambda(h)\in\sigma(\CA_h)\,$. We are also looking for uniform resolvent estimates in any half-plane free of eigenvalues.\\
One of the main difficulties of this task is that, due to possible pseudospectral effects, a quasimode construction may not be sufficient to locate
an eigenvalue.\\

The question considered here is related to stability problems for equations of the form
\beq\label{eqGdDom}
\left\{
\begin{array}{ll}
 \pa_t\psi_R-\Delta\psi_R+iRV(x/R)\psi_R = \lambda_R\psi_R, & (t,x)\in(0,+\infty)\times\Omega_R\,,\\
 \psi_R(t,x) = 0, & (t,x)\in(0,+\infty)\times\pa\Omega_R\,,\\
 \psi_R(0,x) = \psi_R^0(x), & x\in\Omega_R\,,
 \end{array}
 \right.
\eeq
where $\Omega_R=\{Rx : x\in\Omega\}\,$, in the large domain limit $R\rightarrow+\infty\,$.
This system can be interpreted as a linearization of the time-dependent Ginzburg-Landau system in superconductivity,
without magnetic field and in a large smooth domain. From this point of view, the following results should be compared with those of
\cite{Alm,AH,AHP1,AHP2,AHP3}.\\
Similar questions have also been considered in \cite{BHHR} in a $1$-dimensional setting to understand the controllability of some degenerate parabolic
equations.\\
In addition to these applications, the results stated in this paper might have some independent, theoretical interest in the growing field of
non-selfadjoint spectral theory.\\

We shall first focus on the case where the potential $V$ has no critical point. Here again, this assumption makes sense in the framework of
superconductivity, see \cite{Alm} and Section \ref{sSuper}. More precisely, we will prove the following:
\begin{theorem}\label{thmNSAschro1}
Let $n\geq1$ and $V\in\mathcal{C}^\infty(\bar\Omega ; \mathbb{R})$ be such that,
for every $x\in\bar\Omega$, $\nabla V(x)\neq0\,$. Let
\beq\label{defPaPerp}
 \pa\Omega_\perp = \{x\in\pa\Omega : \nabla V(x)\times\vec n(x) = 0\}\,,
\eeq
where $\vec n(x)$ denotes the outward normal on $\pa\Omega$ at $x$\,.
\begin{enumerate}[(i)]
 \item Assume that $\pa\Omega_\perp\neq\emptyset\,$. Let $\mu_1<0$ be the rightmost zero of the Airy function $Ai$\,,
and let
 \beq\label{defJm}
  J_m = \min_{x\in\pa\Omega_\perp}|\nabla V(x)|\,.
 \eeq
Then we have
 \beq\label{limSpect1}
 \varliminf\limits_{h\to0}\frac{1}{h^{2/3}}\inf\Re\sigma(\CA_h) \geq \frac{|\mu_1|}{2}J_m^{2/3}\,,
 \eeq
 where $\CA_h$ is the operator defined by (\ref{defAh})\,.\\
 Moreover, for every $\varepsilon>0\,$, there exists $h_\varepsilon\in(0,h_0)$ and $C_\varepsilon>0$ such that
 \beq\label{estRes1}
 \forall h\in(0,h_\varepsilon),~~~
 \sup_{\tiny{\begin{array}{c}\gamma\leq|\mu_1|J_m^{2/3}/2,\\ \nu\in\mathbb{R}\end{array}}}
 \|(\CA_h-(\gamma -\varepsilon)h^{2/3}-i\nu)^{-1}\|\leq\frac{C_\varepsilon}{h^{2/3}}\,.
 \eeq
 \item Assume that $\pa\Omega_\perp = \emptyset$\,, then 
 \[
\lim_{h\rightarrow0}\frac{1}{h^{2/3}}\inf\Re\sigma(\CA_h) = +\infty\,,  
 \]
and for all $\omega\in\mathbb{R}\,$, there exists $h_\omega>0$ and 
 $C_\omega'>0$ such that
  \beq
 \forall h\in(0,h_\omega),~~~
 \sup_{\tiny{\begin{array}{c}\gamma\leq\omega,\\ \nu\in\mathbb{R}\end{array}}}
 \|(\CA_h-\gamma h^{2/3}-i\nu)^{-1}\|\leq\frac{C_\omega'}{h^{2/3}}\,.
 \eeq
\end{enumerate}
\end{theorem}
This result is essentially a reformulation of those stated in \cite{Alm}, but the proof presented here, based on locally approximating models,
gives a good overview of the underlying phenomena involved and might be more convenient for possible generalizations of this statement.\\
As we shall see in the proof of this first statement, we will not be able to prove that $\frac{|\mu_1|}{2}J_m^{2/3}$ is the exact limit
for $h^{-2/3}\inf\Re\sigma(\CA_h)$ as $h\rightarrow0$\,. This is because we will have to approximate $\CA_h$ in the neighborhood of 
$\pa\Omega_\perp$ by operators whose resolvents are not compact for $n\geq2$. However, this result can still be used to 
obtain some decay estimates for equations of the form (\ref{eqGdDom}), see Corollary \ref{corNSAschroSG} and Sections \ref{sDecaySG} and \ref{sSuper}.\\
In dimension $1$, obviously, this problem of non-compact resolvent will not appear, hence we can state a more accurate result:
\begin{theorem}\label{thmNSAschro1_1d}
Let $h_0>0$, $a,b\in\mathbb{R}$, $a<b$, and $V\in\mathcal{C}^\infty((a,b) ; \mathbb{R})\,$. For $h\in(0,h_0)\,$, let
\[
 \CA_h = -h^2\frac{d^2}{dx^2}+iV(x)\,,~~~\CD(\CA_h) = H_0^1(a,b)\cap H^2(a,b)\,.
\]
Assume that, for every $x\in(a,b)$, $V'(x)\neq0$\,. Then,
 \beq\label{limSpect1d}
 \lim\limits_{h\to0}\frac{1}{h^{2/3}}\inf\Re\sigma(\CA_h) = \frac{|\mu_1|}{2}J^{2/3}\,,
 \eeq
 where  $J = \min(|V'(a)|,|V'(b)|)$ and $\mu_1$ denotes the rightmost zero of the Airy function $Ai$\,.
\end{theorem}
The problem of optimality in (\ref{limSpect1}), in the general, $n$-dimensional setting, is left for future considerations.\\

In the case where the potential $V$ has critical points in $\Omega$, the spectrum of $\CA_h$ is expected to behave
differently. The following statement shows that the quantity $\inf\Re\sigma(\CA_h)$ is no longer
determined by the behavior at the boundary, but by the shape of the potential near
the critical points.

\begin{theorem}\label{thmNSAschro2}
Let $V$ be a Morse function on $\bar\Omega\,$, without critical point in $\pa\Omega$ and with at least one critical point in
$\Omega\,$. Let $x_1^c,\dots,x_p^c\,$, $p\in\mathbb{N}^*\,$, denote those critical points,
and for $k=1,\dots,p$, let
\beq\label{defKappak}
 \kappa_k = \sum_{j=1}^n\sqrt{|\lambda_j^k|}\,, 
\eeq
where $\{\lambda_j^k\}_{j=1,\dots,n} = \sigma(\mathrm{Hess }V(x_k^c))\,$.\\ Let
\[
 \kappa = \min_{k=1,\dots,p}\kappa_k\,,
\]
and assume that, if $\kappa_k = \kappa\,$, then for any $\ell\neq k\,$,
\beq\label{hypNiveaux}
V(x_k^c)\neq V(x_\ell^c)\,.
\eeq
Then,
\beq\label{limSpect2}
\lim\limits_{t\to0}\frac{1}{h}\inf\Re\sigma(\CA_h) = \frac{\kappa}{2}\,.
\eeq
Moreover, for every $\varepsilon>0\,$, there exists $h_\varepsilon\in(0,h_0)$ and $C_\varepsilon>0$ such that
 \beq\label{estRes2}
 \forall h\in(0,h_\varepsilon),~~~
 \sup_{\tiny{\begin{array}{c}\gamma\leq\kappa/2,\\ \nu\in\mathbb{R}\end{array}}}
\|(\CA_h-(\gamma-\varepsilon) h-i\nu)^{-1}\|\leq\frac{C_\varepsilon}{h}\,.
 \eeq
 \end{theorem}
The assumption (\ref{hypNiveaux}) is meant to avoid any resonance phenomenon between two wells. Note that, unlike in Theorem~\ref{thmNSAschro1},
here we give the exact limit for $h^{-1}\inf\Re\sigma(\CA_h)$\,.\\

As mentioned above, the previous theorems enable us to state some decay estimates for the semigroup associated with $\CA_h$.
\begin{corollary}\label{corNSAschroSG}
For all $\varepsilon>0$, there exists $h_\varepsilon\in(0,h_0)$ and $M_\varepsilon>0$ such that:
\begin{enumerate}[(i)]
 \item Under the assumptions of Theorem \ref{thmNSAschro1},
 \beq\label{estSG1}
 \forall h\in(0,h_\varepsilon)\,,~~\forall t>0\,,~~~\|e^{-t\CA_h}\|_{\mathcal{L}(L^2(\Omega))} \leq M_\varepsilon
 \exp(-(|\mu_1|J_m^{2/3}/2-\varepsilon)h^{2/3}t)\,.
 \eeq
 \item Under the assumptions of Theorem \ref{thmNSAschro2},
  \beq\label{estSG2}
 \forall h\in(0,h_\varepsilon)\,,~~\forall t>0\,,~~~\|e^{-t\CA_h}\|_{\mathcal{L}(L^2(\Omega))} \leq M_\varepsilon
 \exp(-(\kappa/2-\varepsilon)ht)\,.
 \eeq
 \item Under the assumptions of Theorem \ref{thmNSAschro1_1d}, the constant $|\mu_1|J_m^{2/3}/2$ is optimal in (\ref{estSG1}),
 as well as the exponent of $h$. Similarly, under the assumptions of Theorem \ref{thmNSAschro2}, the constant $\kappa/2$ is optimal
 in (\ref{estSG2}), as well as the exponent of $h$\,.
\end{enumerate}
\end{corollary}
This corollary will follow easily from Theorems \ref{thmNSAschro1}, \ref{thmNSAschro1_1d} and \ref{thmNSAschro2}, by using a refined, quantitative
version of the Gearhardt-Pr\"uss Theorem, see \cite{HelSjo}.\\

Many interesting questions, which arise naturally in superconductivity theory, are left aside from this paper and should be investigated in
future research.
First of all, as recalled in Section \ref{sSuper}, the time-dependent Ginzburg-Landau equations involve a non-linear term of the form $(1-|\psi|^2)\psi$,
which shall not be considered in this work. 
The recent work of Y. Almog and B. Helffer \cite{AH} includes the analysis of this
non-linearity in the presence of a magnetic field, but as far as we know, this non-linear problem has not been considered yet
in the simpler case where the magnetic field
is neglected.\\
Secondly, here we only consider the case of a smooth domain $\Omega$. As explained in \cite{Alm}, most physically relevant domains would instead
contain some singularities, such as corners with right-angles. However, since Y. Almog \cite{Alm} has already
considered this feature under
the assumption of a potential without critical point, and since the case of a Morse potential is outside the scope of superconductivity theory, 
this question shall not be considered here. Nevertheless, our guess is that the results stated in Theorems~\ref{thmNSAschro1} and
\ref{thmNSAschro2} would be similar
for a domain with right-angled corners at the boundary, and that the proof could be easily adjusted by adding a model acting on a quarter of space in order
to approximate the operator $\CA_h$ near those singularities.\\
Finally, we think that it could be interesting to analyze the effect of a magnetic field when the electric potential has critical points. Namely, the
behavior of $\inf\Re(\CA_{\mathbf{A},V,h})$ should be investigated, where
\[
 \CA_{\mathbf{A},V,h} = -(h\nabla-i\mathbf{A}(x))^2+iV(x),
\]
and where $V$ satisfies the assumptions of Theorem \ref{thmNSAschro2}. This problem has been considered in \cite{AH} in the case where the electric
potential $V$ has no critical point. Of course, some additional conditions on the magnetic field $\mathbf{B} = \curl\mathbf{A}$ should be added 
in order to understand this question.\\ 

Section \ref{sModels} is dedicated to the analysis of some simplified models which shall be used as local approximations for operator $\CA_h$.
In Section \ref{sCoord}, we locally straighten the boundary by introducing a system of local coordinates, previously used in \cite{AH,FouHel,Pan}.
We prove Theorem \ref{thmNSAschro1} in Section \ref{sthm1}, and the lower bound of Theorem \ref{thmNSAschro2} in Section \ref{sthm3Lower}.
We complete the proof of Theorem \ref{thmNSAschro1_1d} (upper bound) in Section \ref{sthm2}, and we prove the upper bound of Theorem 
\ref{thmNSAschro2} in Section \ref{sthm3Upper}.
Section \ref{sDecaySG} is devoted to the proof of Corollary \ref{corNSAschroSG}. Finally, in Section \ref{sSuper},
we give a possible application for the previous results in superconductivity theory, recovering the results of \cite{Alm}.

\section{Simplified models}\label{sModels}
In this section, we consider the simplified cases where $\Omega$ is either the whole space $\mathbb{R}^n$, or the half-space
\beq\label{defHalfSpace}
\mathbb{R}_+^n = \{x=(x_1,\dots,x_n)\in\mathbb{R}^n : x_n>0\}.
\eeq
Furthermore the potential $V$ will be assumed to be a linear function or,
in Subsection \ref{ssModelQuad}, a quadratic form.\\
In Sections \ref{sthm1} to \ref{sthm3Upper}, we shall use these simplified models as local approximations of the more general operator $\CA_h\,$.

\subsection{Whole space model, and particular half-space models}\label{ssModels}
In this subsection, we mainly refer to \cite{Hel2}, and reformulate the $2$-dimensional statements therein in the $n$-dimensional setting.\\
We shall consider three model operators $-\Delta+i\ell\,$, where $\ell(x)=J\cdot x$ is a linear function: the first one in $\mathbb{R}^n\,$, the second one in
$\mathbb{R}_+^n$ with $J$ parallel to $\pa\mathbb{R}_+^n\,$, and the third one in $\mathbb{R}_+^n$ with $J$ orthogonal to $\pa\mathbb{R}_+^n\,$.

\subsubsection{The whole space model}\label{pWhole}
Let $J=(J_1,\dots,J_n)\in\mathbb{R}^n$ and $$\CA_0 = -\Delta+i\ell $$ acting on $L^2(\mathbb{R}^n)\,$,
where $\ell(x) = J\cdot x\,$. Up to an orthogonal change of variable followed
by the scale change $x\mapsto|J|^{1/3}x$, we can assume that $\CA_0$ has the form $\CA_0 = -\Delta+ix_1\,$.\\
Recalling that the complex Airy operator $-\frac{d^2}{dx_1^2}+ix_1$ on $L^2(\mathbb{R})$ has empty spectrum, we then get as in \cite{Hel2},
Proposition $7.1$,
\begin{lemma}\label{lemModelRn}
We have $\sigma(\CA_0) = \emptyset$, and for all $\omega\in\mathbb{R}$, there exists $C_\omega^0$ such that
\beq
\sup_{\Re z\leq\omega}\|(\CA_0-z)^{-1}\|\leq C_\omega^0.
\eeq
\end{lemma}

\subsubsection{Parallel current in the half-space}\label{pParallel}

Now we consider the Dirichlet (resp. Neumann) realization $\CA_\sslash^D$ (resp. $\CA_\sslash^N$)
of $-\Delta+i(J_1x_1+\dots+J_{n-1}x_{n-1})$ in $\mathbb{R}_+^n\,$. 
As in \cite{Hel2}, Subsection $7.3$, we can use the decomposition $L^2(\mathbb{R}^n)=\mathfrak{P}\oplus\mathfrak{I}$, where $\mathfrak{P}$ and
$\mathfrak{I}$ denote respectively the even and odd functions in $L^2(\mathbb{R}^n)$ with respect to the $x_n$ variable, and check that
$\sigma(\CA_0) = \sigma(\CA_\sslash^D)\cup\sigma(\CA_\sslash^N)$. Hence in view of Lemma \ref{lemModelRn}, the spectra of $\CA_\sslash^D$ and
$\CA_\sslash^N$ are empty.\\
Moreover, since $\CA_\sslash^D$ (resp. $\CA_\sslash^N$) is the restriction of $\CA_0$ to $\mathfrak{P}$ (resp. $\mathfrak{I}$), the
resolvent estimate in Lemma \ref{lemModelRn} yields
\begin{lemma}\label{lemModelParallel}
$\sigma(\CA_\sslash^D) = \sigma(\CA_\sslash^N)=\emptyset$, and for all $\omega\in\mathbb{R}$ and $\sharp=D,N$,
\beq
\sup_{\Re z\leq\omega}\|(\CA_\sslash^\sharp-z)^{-1}\|\leq C_\omega^0.
\eeq
\end{lemma}

\subsubsection{Perpendicular current in the half-space}\label{pPerp}
Let $J_n\in\mathbb{R}$ and $\CA_\perp$ be the Dirichlet realization of $-\Delta+iJ_nx_n$ in $L^2(\mathbb{R}_+^n)$. Here again, we can easily
adapt the results of \cite{Hel2} (see Proposition $7.2$ therein) to any dimension $n\geq2$ to obtain
\[
 \sigma(\CA_\perp) = \bigcup_{\lambda\in\sigma(\CA_{x_n})}\{\lambda+r : r>0\},
\]
where $\CA_{x_n}$ denotes the Dirichlet realization of the complex Airy operator $-\frac{d^2}{dx_n^2}+iJ_nx_n$ in $L^2(\mathbb{R}^+)$.
Recalling from \cite{Alm} that $\sigma(\CA_{x_n}) = \{e^{i\pi/3}|\mu_j|J_n^{2/3}\}_{j\geq1}$ where $\mu_j<0$ are the zeroes
of the Airy function $Ai\,$, we get
\begin{lemma}\label{lemModelPerp}
 $\sigma(\CA_\perp) = \{|\mu_j|J_n^{2/3}e^{i\pi/3}+r : j\geq1,~r>0\}$, and for all $\omega<|\mu_1|J_n^{2/3}/2$, there exists $C_\omega^\perp$ such that
 \beq
 \sup_{\Re z\leq\omega}\|(\CA_\perp-z)^{-1}\|\leq C_\omega^\perp.
 \eeq
\end{lemma}

In the following subsection, we consider a case which was not studied in \cite{Hel2}: the operator $-\Delta+iJ\cdot x$ in the half-space, where
neither $J_n$ nor $(J_1,\dots,J_{n-1})$ vanishes.

\subsection{General current in the half-space}\label{ssModelQcq}

Let $J=(J_1,\dots,J_n)\in\mathbb{R}^n$ be such that $J_n\neq0$ and $(J_1,\dots,J_{n-1})\neq0\,$.
We want to study the spectrum and the resolvent of an operator acting on $L^2(\mathbb{R}_+^n)$ as
$-\Delta+iJ\cdot x\,$, with a domain which includes the Dirichlet condition at $x_n=0\,$. The imaginary part
\[
 \ell(x,y) = J\cdot x
\]
of the potential does not have a constant sign, hence we are unable to use the variational approach to define the operator.
We shall instead define the operator by separation of variables.\\

Let $x'=(x_1,\dots,x_{n-1})$ denote the $(n-1)$ first coordinates of a vector $x\in\mathbb{R}^n$. Let
\beq\label{defCA_x'}
 \CA_{x'} = -\Delta_{x'}+iJ'\cdot x'\,,
\eeq
and let $\CA_{x_n}^+$ be the Dirichlet realization in $\mathbb{R}^+$ of the complex Airy operator
\beq\label{defCA_xn}
 -\frac{d^2}{dx_n^2}+iJ_nx_n\,.
\eeq
Both $\CA_{x'}$ and $\CA_{x_n}^+$ are maximal accretive, hence they are generators of contraction semigroups $(e^{-t\CA_{x'}})_{t>0}$ and
$(e^{-t\CA_{x_n}^+})_{t>0}$ respectively. One can easily check that the family $(e^{-t\CA_{x'}}\otimes e^{-t\CA_{x_n}^+})_{t>0}$ is a contraction
semigroup in $L^2(\mathbb{R}_+^n)\,$. Thus, we can define the desired operator as follows:
\begin{definition}
$\CA_+$ is the generator of the semigroup $(e^{-t\CA_{x'}}\otimes e^{-t\CA_{x_n}^+})_{t>0}\,$.
\end{definition}

In order to describe the domain of the operator $\CA_+$, we recall the following statement (see \cite{ReSim}, Theorem X.$49$):
\begin{theorem}\label{thmcoeur}
Let $\CA$ be the generator of a contraction semigroup on a Hilbert space $\CH$. Let $\CD\subset\CD(\CA)$ be a dense subset of $\CH$,
such that $e^{-t\CA}\CD\subset\CD\,$. Then, $\CD$ is a \emph{core} for $\CA\,$, that is
 \[
  \CA = \barr{\CA_{|\CD}}\,.
 \]
\end{theorem}

Let $\CD=\CD(\CA_{x'})\odot\CD(\CA_{x_n}^+)$ be the set of all finite linear combinations of functions of the form $f\otimes g=f(x')g(x_n)$, where 
$f\in\CD(\CA_{x'})$ and $g\in\CD(\CA_{x_n}^+)$. Then it is clear that $\CD$ satisfies the conditions of Theorem \ref{thmcoeur} with
$\CA=\CA_+\,$, hence
$\CA_+ = \barr{{\CA_+}_{|\CD}}\,$. Consequently, we have the following characterization of the domain:
\begin{eqnarray}
 \CD(\CA_+) &=& \{u\in L^2(\mathbb{R}_+^2) : \exists (u_j)_{j\geq1}\in\CD^\mathbb{N}\,,~
 u_j\underset{\tiny{j\rightarrow+\infty}}{\overset{L^2}{\longrightarrow}}u\,,\nonumber\\
 & &\quad\quad \quad\quad  (\CA_+u_j)_{j\geq1}~\textrm{is a Cauchy sequence}~\}\,.\label{caracDomA+}
\end{eqnarray}

We already mentioned that the sequilinear form associated with $\CA_+$ is not coercive. We can indeed consider the sequence
\[
u_j(x',x_n) := \sqrt{\chi(|(x_2,\dots,x_n)|-2)\left(\chi(x_1-j)+\chi(x_1+j)\right)},
\]
where $\chi\in\mathcal{C}_0^\infty(]-2,2[ ; [0,1])$ is equal to $1$ on $[-1,1]$, and check that
$\sca{\CA_+u_j}{u_j}$ is bounded, whereas
 $\||\ell|^{1/2}u_j\|_{L^2(\mathbb{R}_+^2)}\rightarrow+\infty$ as $j\rightarrow+\infty$.\\
 Similarly, we can check that the resolvent of $\CA_+$ is not compact, by considering, for some
$u_0\in~\mathcal{C}_0^\infty(\mathbb{R}_+^n)$, the sequence 
\[
 u_j(x,y):=u_0(x_1-2jJ_2M,x_2+2jJ_1M,x_3,\dots,x_n),
\]
where $M$ denotes the diameter of $\Supp u_0$.
Then the supports of $u_j$ are disjoints and translated along a direction leaving $J\cdot x$ constant. Hence, 
it is straightforward to check that the sequence $(\CA_+u_j)_{j\geq1}$ is bounded while $(u_j)_{j\geq1}$ has no converging subsequence.\\
Finally, one can show by a similar construction that, for $j=1,\dots,n$, we can not expect to have a control of
$\|(-\frac{d^2}{dx_j^2}+iJ_jx_j)u\|_{L^2(\mathbb{R}_+^n)}$ by the graph norm of $\CA_+$.\\
However, we prove in the following lemma that we can control separately $\|u\|_{H^2(\mathbb{R}_+^n)}$ and $\|\ell u\|_{L^2(\mathbb{R}_+^n)}\,$,
which gives a good description of the domain:
\begin{lemma}
We have
\beq\label{descrDomA+}
\CD(\CA_+) = H_0^1(\mathbb{R}_+^n)\cap H^2(\mathbb{R}_+^n)\cap L^2(\mathbb{R}_+^n ; |\ell(x)|^2dx)\,,
\eeq
and there exists $C>0$ such that, for all $u\in\CD(\CA_+)\,$,
\beq\label{estDomA+}
\|\Delta u\|_{L^2(\mathbb{R}_+^n)}^2 + \|\ell u\|_{L^2(\mathbb{R}_+^n)}^2 \leq \|\CA_+ u\|_{L^2(\mathbb{R}_+^n)}^2+
C\|\nabla u\|_{L^2(\mathbb{R}_+^n)}\|u\|_{L^2(\mathbb{R}_+^n)}\,.
\eeq
\end{lemma}
\textbf{Proof: }
We use the characterization (\ref{caracDomA+}). Let $u\in\CD(\CA_+)$ and
$(u_j)_{j\geq1}\in\CD^\mathbb{N}$ such that $u_j\underset{\tiny{j\rightarrow+\infty}}{\overset{L^2}{\longrightarrow}}u$
and $(\CA_+u_j)_{j\geq1}$ is a Cauchy sequence. Then, using the identity
\[
 \Re\sca{\CA_+v}{v} = \|\nabla v\|_{L^2(\mathbb{R}_+^n)}^2\,,
\]
we see that $(\nabla u_j)_{j\geq1}$ is a Cauchy sequence in $L^2(\mathbb{R}_+^n)$, hence
\beq\label{cvNabla}
u_j\underset{\tiny{j\rightarrow+\infty}}{\overset{H^1}{\longrightarrow}} u\,,
\eeq
and $u\in H_0^1(\mathbb{R}_+^n)\,$.\\
In order to prove (\ref{estDomA+}), we write (all the norms denoting $L^2$ norms)
\begin{eqnarray}
 \|\CA_+u_j\|^2 & = & \sca{(-\Delta+i\ell)u_j}{(-\Delta+i\ell)u_j} \nonumber\\
 & = & \|\Delta u_j\|^2+ \|\ell u_j\|^2 + 2\Im\sca{-\Delta u_j}{\ell u_j}\,.\label{estDomInter}
\end{eqnarray}
Besides, we have
\begin{eqnarray*}
\Im \sca{-\Delta u_j}{\ell u_j} & = & \Im\int_{\mathbb{R}_+^n}\nabla u_j(x)\cdot\barr{\nabla (\ell u_j)(x)}dx\\
 & = & \Im\left(\int_{\mathbb{R}_+^n}\ell(x)|\nabla u_j(x)|^2dx + \int_{\mathbb{R}_+^n}\nabla u_j(x)\cdot\barr{\nabla\ell(x)}\barr{u_j(x)}dx\right) \\
 & = & \Im\int_{\mathbb{R}_+^n} J\cdot\nabla u_j(x)\barr{u_j(x)}dx\,.
\end{eqnarray*}
Hence, for some $C>0\,$,
\[
 |\Im \sca{-\Delta u_j}{\ell u_j}|\leq C\,\|\nabla u_j\|\, \|u_j\|\,.
\]
Thus, according to (\ref{estDomInter}), estimate (\ref{estDomA+}) holds for the functions $u_j$. Consequently,
$(u_j)_{j\geq1}$ is a Cauchy sequence in $H^2(\mathbb{R}_+^n)$ and in $L^2(\mathbb{R}_+^n ; |\ell(x)|^2dx)$, and (\ref{descrDomA+}) follows,
as well as (\ref{estDomA+}) for every $u\in\CD(\CA_+)\,$.
\hfill $\boxminus$\\

Now we answer the question of the spectrum of $\CA_+\,$. Since $\CA_{x'}$ has empty spectrum (see Subsection \ref{ssModels}), we expect
$\sigma(\CA_+)$ to be empty as well. In order to prove it, we use semigroup estimates.
\begin{proposition}\label{propSpecA+}
We have $\sigma(\CA_+) = \emptyset\,$. Moreover, for every $\omega\in\mathbb{R}\,$, there exists $C_\omega>0$ such that
\beq\label{ResA+}
\sup_{\tiny{\Re z \leq \omega}}\|(\CA_+-z)^{-1}\| \leq C_\omega\,.
\eeq
Finally, the semigroup generated by $\CA_+$ satisfies
\beq\label{decSGA+}
\forall t>0,~~\|e^{-t\CA_+}\|\leq e^{-(n-1)t^3/12}\,.
\eeq
\end{proposition}

\textbf{Proof: }
Let us recall that $e^{-t\CA_+} = e^{-t\CA_{x'}}\otimes e^{-t\CA_{x_n}^+}\,$, where $\CA_{x'}$ and $\CA_{x_n}^+$ are respectively
defined by (\ref{defCA_x'}) and (\ref{defCA_xn}). We also recall the following estimates (see \cite{Hel2,Mart}):
\beq\label{estSGAiry}
 \forall t>0\,,~~~ \|e^{-t\CA_{x'}}\| = e^{-(n-1)t^{3}/12}\,,
\eeq
and for all $\omega<|\mu_1|/2\,$, where $\mu_1$ is the rightmost zero of the Airy function, there exists $M_\omega>0$ such that
\beq\label{estSGAiry+}
 \forall t>0\,,~~~ \|e^{-t\CA_{x_n}^+}\|\leq M_\omega\, e^{-\omega t}\,.
\eeq
Thus, (\ref{decSGA+}) follows, and the formula
\beq\label{formResSG}
 (\CA_+-z)^{-1} = \int_0^{+\infty}e^{-t(\CA_+-z)}dt\,,
\eeq
which holds \emph{a priori} for $\Re z<0\,$, can be extended to the whole complex plane.
Hence the resolvent of
$\CA_+$ is an entire function, and we have $\sigma(\CA_+) = \emptyset$ as well as (\ref{ResA+})\,.
\hfill $\square$

\subsection{Uniform resolvent estimate with respect to the angle}
In the proof of the main theorems, we will need to manage the transition between the case of an orthogonal current and a general
transverse current. Let $J_0\,, J_1>0$ such that $J_0<J_1\,$. Then, for $J\in\mathbb{R}$ such that $J_0<|J|<J_1\,$, $\theta\in[0,\pi]$ and
$\vec v\in \mathbb{R}^{n-1}\,$, $|\vec v|=1\,$, we set
\beq\label{defAJtheta}
\CA_+(J,\theta,\vec v) = -\Delta+iJ(\sin\theta~\vec v\cdot x'+\cos\theta~x_n)\,,
\eeq
acting on $L^2(\mathbb{R}_+^n)\,$.
Then, using the results of previous subsections, we shall prove that:
\begin{lemma}\label{lemAngle}~
 \begin{enumerate}[(i)]
  \item For all $\varepsilon>0\,$, there exists $K_\varepsilon>0$ such that, 
  for all $J$ satisfying $J_0\leq|J|\leq J_1$  and $\vec v\in \mathbb{R}^{n-1}$ such that $|\vec v|=1\,$,
  \beq\label{resAngle1}
  \sup_{\tiny{\begin{array}{c}\theta\in[0,\pi]\\ \Re z\leq|\mu_1|/2 \end{array}}}\|(\CA_+(J,\theta,\vec v)-(z-\varepsilon)|J|^{2/3})^{-1}\|\leq
  \frac{K_\varepsilon}{|J|^{2/3}}\,.
  \eeq
  \item Let $\theta_0\in(0,\pi/2)\,$. Then,
  for all $\omega\in\mathbb{R}\,$, there exists $K_\omega'>0$ such that, for every $J$ satisfying
  $J_0\leq|J|\leq J_1$ and $\vec v\in \mathbb{R}^{n-1}\,$, $|\vec v|=1\,$,
  \beq\label{resAngle2}
  \sup_{\tiny{\begin{array}{c}\theta_0\leq\theta\leq\pi-\theta_0
  \\ \Re z\leq\omega \end{array}}}\|(\CA_+(J,\theta,\vec v)-z|J|^{2/3})^{-1}\|\leq
  \frac{K_\omega'}{|J|^{2/3}}\,.
  \eeq
 \end{enumerate}
\end{lemma}
\textbf{Proof: }
For $j=1,\dots,n-1\,$, let
\[
 \CA_j(J,\theta,v_j) = -\frac{d^2}{dx_j^2} + iJv_j\sin\theta~x_j\,,
\]
acting on $L^2(\mathbb{R})\,$, and 
\[
 \CA_n^+(J,\theta) = -\frac{d^2}{dx_n^2}+iJ\cos\theta~x_n\,,
\]
acting on $L^2(\mathbb{R}^+)\,$.\\
 Using that each $\CA_j(J,\theta,v_j)\,$, $j=1,\dots,n-1$ and $\CA_n^+(J,\theta)$ is maximal accretive, we can easily
check that, for $\vec v = (v_1,\dots,v_{n-1})\,$,
\beq\label{OtimesAJthetav}
 e^{-t\CA(J,\theta,\vec v)} = e^{-t\CA_1(J,\theta,v_1)}\otimes\dots\otimes e^{-t\CA_{n-1}(J,\theta,v_{n-1})}\otimes e^{-t\CA_n^+(J,\theta)}\,.
\eeq
Since $|\vec v|=1\,$, we can choose $k$ such that $|v_k|\geq1/\sqrt{n}\,$. We can also assume that $Jv_k\sin\theta$ and
$J\cos\theta$ are both non-negative (if not, replace $\CA_k(J,\theta,v_k)$ or $\CA_n^+(J,\theta)$ by its adjoint).\\
Notice now that, for $\theta\in(0,\pi)\,$, by rescaling $x\mapsto (Jv_k\sin\theta)^{1/3}x\,$, we have
\[e^{-t\CA_k(J,\theta,v_k)} = e^{-t(|J|v_k\sin\theta)^{2/3}\CA_k(1,\pi/2,1)}\,,\] and by rescaling $x\mapsto (|J|\cos\theta)^{1/3}x$, we have similarly
\[e^{-t\CA_n^+(J,\theta)} = e^{-t(|J|\cos\theta)^{2/3}\CA_n^+(1,0)}\,.\]
Hence, if $\theta\in[\theta_0,\pi-\theta_0]\,$, then according to (\ref{estSGAiry}) and (\ref{OtimesAJthetav}),
\[
 \|e^{-t\CA(J,\theta,\vec v)}\| \leq \|e^{-t(|J|v_k\sin\theta)^{2/3}\CA_k(1,\pi/2,1)}\| = e^{-t^3|Jv_k\sin\theta|^2/12}
 \leq e^{-\varepsilon_0^2t^3/12n}\,,
\]
where $\varepsilon_0 = \sin\theta_0\,$.\\
 Thus, formula (\ref{formResSG}) yields (\ref{resAngle2}).\\
In order to prove (\ref{resAngle1}), for $\omega<|\mu_1|/2\,$, we write, using (\ref{estSGAiry}), (\ref{estSGAiry+}) and (\ref{OtimesAJthetav}),
\begin{eqnarray*}
 \|e^{-t\CA(J,\theta,\vec v)}\| & \leq & \|e^{-t(|J|v_k\sin\theta)^{2/3}\CA_k(1,\pi/2,1)}\|\|e^{-t(|J|\cos\theta)^{2/3}\CA_n^+(1,0)}\|,\\
  & \leq & M_\omega\exp\left(-\frac{t^3}{12n}|J|^2\sin^2\theta-t|J|^{2/3}\omega(\cos\theta)^{2/3}\right) \\
   & \leq & M_\omega e^{-|J|^{2/3}\omega t}e^{-g_\theta(t|J|^{2/3})}\,,
\end{eqnarray*}
where
\[
 g_\theta(s) = -\frac{1}{12n}\sin^2\theta s^3+\omega s(1-(\cos\theta)^{2/3})\,.
\]
It is then straightforward to check that $g_\theta$ is bounded in $\mathbb{R}^+\,$, uniformly with respect to $\theta\in(0,\pi)\,$. Hence, (\ref{resAngle1})
follows from formula (\ref{formResSG}).
\hfill $\boxminus$
\subsection{Quadratic potential in the whole space}\label{ssModelQuad}
In order to prove Theorem \ref{thmNSAschro2}, we will need to understand the pseudospectral behavior of operators of the form
$-\Delta+i\mathcal{Q}$ acting in $L^2(\mathbb{R}^n)$, where $\mathcal{Q}$ is a quadratic form.\\
More precisely, let $\lambda = (\lambda_1,\dots,\lambda_n)\in\mathbb{R}^n$ such that $\lambda_j\neq0\,$, $j=1,\dots,n\,$, and
\beq\label{defQuad}
 \CH_\mathcal{Q} = -\Delta+i\mathcal{Q}\,,~~~\mathcal{Q}(x)=\mathcal{Q}_\lambda(x) = \sum_{j=1}^n\lambda_jx_j^2\,.
\eeq
We want to determine the spectrum of $\CH_\mathcal{Q}$ and to control its resolvent uniformly on any half-plane included in the resolvent set.\\

For any $\alpha\in\mathbb{R}\setminus\{0\}\,$, let
\[
 \CH_\alpha = -\frac{d^2}{dx^2}+i\alpha x^2,~~~\CD(\CH_\alpha) = H^2(\mathbb{R})\cap L^2(\mathbb{R} ; x^4dx)\,,
\]
be the ($1$-dimensional) complex harmonic oscillator \cite{Boul,Dav2,DavKui,PraSta1}. Let us recall that
\beq\label{specDavies}
 \sigma(\CH_\alpha) = \{(2k+1)\sqrt{|\alpha|}e^{\pm i\pi/4} : k\in\mathbb{N}\},
\eeq
where $\pm = \mathrm{sign }~\alpha$\,.\\
 If $\alpha<0\,$, notice indeed that $\CH_\alpha = \CH_{|\alpha|}^*\,$, hence 
$\lambda\in\sigma(\CH_\alpha)$ if and only if $\bar\lambda\in\sigma(\CH_{|\alpha|})\,$. Moreover, for every $\omega<\sqrt{|\alpha|/2}\,$,
there exists $c_\omega>0$ such that
\beq\label{resDavies}
\sup_{\Re z\leq\omega}\|(\CH_\alpha-z)^{-1}\|\leq c_\omega\,,
\eeq
see \cite{Hel1}, Proposition $14.13$, and \cite{Boul,PraSta1}.\\
Now, notice that
\[
 \CH_\mathcal{Q} = \barr{\sum_{j=1}^n I\otimes\dots\otimes I\otimes\CH_{\lambda_j}\otimes I\otimes\dots\otimes I}\,,
\]
(use for instance Theorem \ref{thmcoeur} to check that the domains coincide).\\
Unlike in Subsection \ref{ssModelQcq}, separation of variables is very efficient for $\CH_\mathcal{Q}$ because
the operators $\CH_{\lambda_j}$ appearing in its decomposition are sectorial.
We can indeed apply the spectral mapping theorem due to Ichinose, given in \cite{ReSim}, XIII. $9$, which
yields
\beq\label{SMTquadratic}
\sigma(\CH_\mathcal{Q}) = \sigma(\CH_{\lambda_1})+\dots +\sigma(\CH_{\lambda_n})\,.
\eeq
In view of (\ref{specDavies}) and (\ref{resDavies}), we then get
\begin{lemma}\label{lemModelQuad}
Let $\lambda = (\lambda_1,\dots,\lambda_n)\in\mathbb{R}^n$ such that $\lambda_j\neq0$ for all $=1,\dots,n\,$, and let
$\sigma_j=\mathrm{sign }~\lambda_j\,$. Then,
\beq
\sigma(\CH_\mathcal{Q}) = \left\{\sum_{j=1}^n(2k_j+1)\sqrt{|\lambda_j|}e^{i\sigma_j\pi/4} : (k_1,\dots,k_n)\in\mathbb{N}^n\right\}\,.
\eeq
Moreover, for all $\omega<\sqrt{|\lambda_1|/2}+\dots+\sqrt{|\lambda_n|/2}\,$, there exists $K_\omega>0$ such that
\beq\label{estResModelQuad}
\sup_{\Re z\leq\omega}\|(\CH_\mathcal{Q}-z)^{-1}\|\leq K_\omega\,.
\eeq
\end{lemma}

\section{Local coordinates near the boundary}\label{sCoord}
In this section, we introduce local coordinates in the neighborhood of some point $b\in\pa\Omega$,
in order to straighten a portion of the boundary. These coordinates will allow us to use the models of previous section as approximate operators
for $\CA_h$.\\
Throughout this section, we mainly refer to \cite{FouHel}, appendix F and \cite{Pan}, although these coordinates have also been used in \cite{AH} in a
$2$-dimensional setting.\\

Let $b\in\pa\Omega$ be fixed. Then, for some neighborhood $\omega\subset\Omega$ of $b$
and some neighborhood $\mathcal{U}$ of the origin in $\mathbb{R}^{n-1}$, there exists a diffeomorphism
\[
 \varphi : \begin{array}{ccc} \mathcal{U} & \longrightarrow & \pa\Omega\cap\omega \\
            y=(y_1,\dots,y_{n-1}) & \longmapsto & x=\varphi(y)
           \end{array}
\]
with $\varphi(0)=b\,$.\\
Then, in these coordinates, the metric induced on $\pa\Omega$ by the euclidian metric of $\mathbb{R}^n$
writes
\[
\sum_{i,j = 1}^{n-1}g_{ij}(y)dy_i\otimes dy_j\,,
\]
where
\beq\label{defg}
\forall i,j=1,\dots, n-1,~~~g_{ij}(y) = \pa_{y_i}\varphi(y)\cdot\pa_{y_j}\varphi(y)\,.
\eeq
We can choose $\varphi$ so that, for every $i,j\in\{1,\dots,n-1\}\,$,
$\pa_{y_i}\varphi(0)\cdot\pa_{y_j}\varphi(0) = \delta_{i,j}\,.$\\
Hence, if $g$ denotes the matrix $g = (g_{ij})_{i,j}$, then $\det g(0) = 1$\,.\\

Now we define some local coordinates in a neighborhood of $b$ in $\Omega\,$.\\
Let $\vec\nu(y)=-\vec n(\varphi(y))\,$, where $\vec n(x)$ is the outward normal of $\pa\Omega$ at $x\,$. We then define the map $\mathcal{F}$ by
\beq\label{defF}
\mathcal{F}(y,z) = \varphi(y)+z\vec\nu(y)\,.
\eeq
Notice that $z = d(\mathcal{F}(y,z),\pa\Omega)\,$.\\
After taking possibly a smaller $\omega$, there exists $z_0>0$ such that
$\mathcal{F}$ is a diffeomorphism from $\mathcal{U}\times(0,z_0)\subset\mathbb{R}^n$
onto $\Omega\cap\omega\,$.\\
In the following we use the notation $\pa_j = \pa_{y_j}$ if $j\in\{1,\dots,n-1\}\,$, and $\pa_n = \pa_z\,$. In the coordinates $(y,z)$ the euclidian
metric in $\mathbb{R}^n$ writes
\[
 \sum_{k=1}^ndx_k\otimes dx_k = \sum_{i,j=1}^{n-1} G_{ij}(y,z)dy_i\otimes dy_j + dz\otimes dz\,,
\]
where for every $i,j=1,\dots,n\,$,
\[
 G_{ij}(y,z) = \pa_i\mathcal{F}(y,z)\cdot\pa_j\mathcal{F}(y,z)\,.
\]
Indeed, we have, for all $i=1,\dots,n-1\,$,
\[
 G_{in}(y,z) = (\pa_i\varphi(y)+z\pa_i\vec\nu(y))\cdot\vec\nu(y) = 0\,,
\]
since $\pa_i\varphi(y)$ and $\pa_i\vec\nu(y)$ are tangent to $\pa\Omega\,$, and
\[
 G_{nn}(y,z) = |\vec\nu(y)|^2 = 1\,.
\]
Besides, for all $i,j\in\{1,\dots,n-1\}\,$, we have
\begin{eqnarray*}
  G_{ij}(y,z) & = & (\pa_i\varphi(y)+z\pa_i\vec\nu(y))\cdot(\pa_j\varphi(y)+z\pa_j\vec\nu(y)) \\
  & = & g_{ij}(y)+z(\pa_i\vec\nu(y)\cdot\pa_j\varphi(y)+\pa_i\varphi(y)\cdot\pa_j\vec\nu(y)) + z^2\pa_i\vec\nu(y)\cdot\pa_j\vec\nu(y) \\
  & = & g_{ij}(y)-2\Omega_{ij}(y)z+z^2\pa_i\vec\nu(y)\cdot\pa_j\vec\nu(y)\,,
\end{eqnarray*}
where $\Omega_{ij} = -\pa_i\varphi\cdot\pa_j\vec\nu$ are the coefficients of the second fundamental form of $\pa\Omega\,$.\\
Notice that, since we have assumed $g_{ij}(0)=\delta_{ij}$, we have also
\beq\label{G0}
G_{ij}(0,0) = \delta_{ij}~~~~\textrm{ and }~~~~\det G(0,0) = 1\,,
\eeq
where $G=(G_{ij})_{i,j}\,$.\\

We conclude this section by giving the expression of operator $\CA_h$ in coordinates $(y,z)\,$. If we set
\beq\label{transfoF}
T_{\mathcal{F}} : \begin{array}{ccc} L^2(\Omega\cap\omega) & \longrightarrow & L^2(\mathcal{U}\times(0,z_0)) \\
            u & \longmapsto & u\circ\mathcal{F}
           \end{array}\,,
\eeq
then
\begin{eqnarray}
T_{\mathcal{F}}(-h^2\Delta_x+iV(x))T_{\mathcal{F}}^{-1} & =& -h^2\sum_{i,j=1}^n\frac{1}{\sqrt{\det G(y,z)}}\pa_i\left(\sqrt{\det G(y,z)} G^{ij}(y,z)
\pa_j\right)\nonumber\\
&& +iV\circ\mathcal{F}(y,z)\,,
 \label{DeltaG}
\end{eqnarray}
where $(G^{ij})_{ij} = (G_{ij})_{ij}^{-1} = G^{-1}\,$. \\
Notice that, since
\[
 G(y,z) = \left(\begin{array}{cccc}
                 & & & 0 \\
        & (G_{ij}(y,z))_{i,j\leq n-1} & & \vdots \\
        & & & 0 \\
        0 & \dots & 0 & 1
                \end{array}
\right)\,,
\]
then we have
\[
 G(y,z)^{-1}  = \left(\begin{array}{cccc}
                 & & & 0 \\
        & (G^{ij}(y,z))_{i,j\leq n-1} & & \vdots \\
        & & & 0 \\
        0 & \dots & 0 & 1
                \end{array}
\right)\,.
\]
Hence, (\ref{DeltaG}) can be reformulated as
\begin{eqnarray}
 T_{\mathcal{F}}(-h^2\Delta_x+iV(x))T_{\mathcal{F}}^{-1} & = &-h^2\left(\sum_{i,j=1}^{n-1}G^{ij}(y,z)\pa_{y_i}\pa_{y_j} +\pa_z^2 \right)\nonumber\\
  & & -h^2\sum_{j=1}^n\beta_j(y,z)\pa_j +V\circ\mathcal{F}(y,z)\,, \label{AhCoordLocales}
\end{eqnarray}
where for all $j\in\{1,\dots,n-1\}$,
\beq\label{defBetaj}
  \beta_j(y,z) = \frac{1}{\sqrt{\det G}}\sum_{i=1}^{n-1}\pa_{y_i}(G^{ij}\sqrt{\det G})~~~~\textrm{ and }~~~~
\beta_n(y,z) = \frac{\pa_z(\sqrt{\det G})}{\sqrt{\det G}}\,.
\eeq
Finally, notice that according to (\ref{G0}), we have
\beq\label{DLGij}
\forall i,j\in\{1,\dots,n-1\}\,,~~~G^{ij}(y,z) = \delta_{ij}+\mathcal{O}(|(y,z)|)\,,~~~|(y,z)|\rightarrow0\,.
\eeq

\section{Lower bound for a potential without critical point}\label{sthm1}
\sectionmark{No critical points: lower bound}
In this section, we work under the assumptions of Theorem \ref{thmNSAschro1}. We prove the results of $(i)$. If $\pa\Omega_\perp = \emptyset\,$, $(ii)$
can be proved alike, by dropping all the terms corresponding to a point $b_j^\perp(h)\in\pa\Omega_\perp$ in the following proof.\\
For $x_0\in\mathbb{R}^n$ and $r>0$, we denote by $B(x_0,\delta)$ the open ball of radius $r$ centered at $x_0\,$. Let
\[
 \pa\Omega_\sslash = \{x\in\pa\Omega : \nabla V(x)\cdot\vec n(x) = 0\}\,,
\]
and
\[
 \pa\Omega_\angle = \pa\Omega\setminus(\pa\Omega_\perp\sqcup\pa\Omega_\sslash)\,.
\]
Our strategy will be to partition the domain $\Omega$ into small subdomains on which $\CA_h$ will be approximated by simpler models based on the operators
studied in Section \ref{sModels}. For some $\rho>0$ to be determined in the following, and for every $h\in(0,h_0)$, we choose
two sets of indices $J_{int}(h)\subset\mathbb{N}\,$, $J_{bdry}(h)\subset\mathbb{N}\,$, and a set of points
\[
 \big\{a_j(h)\in\Omega : j\in J_{int}(h)\big\}\cup\big\{b_k(h)\in\pa\Omega : k\in J_{bdry}(h)\big\}\,,
\]
such that
\[
 \bar\Omega\subset\bigcup_{j\in J_{int}(h)}B(a_j(h),h^\rho)~\cup\bigcup_{k\in J_{bdry}(h)}B(b_k(h),h^\rho)\,,
\]
and such that the closed balls $\bar B(a_j(h),h^\rho/2)\,$, $\bar B(b_k(h),h^\rho/2)$ are all disjoints.\\
Notice that $\sharp J_{int}(h)\propto h^{-n\rho}$ and $\sharp J_{bdry}(h)\propto h^{-(n-1)\rho}$\,.\\
For $\natural = \perp,\sslash,\angle\,$, we define
\[
 J_\natural(h) = \{j\in J_{bdry}(h) : b_j(h)\in\pa\Omega_\natural\}\,.
\]
Now we take a partition of unity in $\Omega\,$,
\[
\left((\chi_{j,h})_{j\in J_{int}(h)},(\zeta_{j,h}^\perp)_{j\in J_\perp(h)},(\zeta_{j,h}^\sslash)_{j\in J_\sslash(h)},
(\zeta_{j,h}^\angle)_{j\in J_\angle(h)}\right)\,,
\]
such that, for every $x\in\bar\Omega\,$,
\beq\label{sommePartition1}
\sum_{j\in J_{int}(h)}\chi_{j,h}(x)^2+\sum_{\natural,~k\in J_\natural(h)}\zeta_{k,h}^\natural(x)^2=1\,,
\eeq
and such that $\Supp \chi_{j,h}\subset B(a_j(h),h^\rho)$ for $j\in J_{int}(h)$, $\Supp\zeta_{j,h}^\natural\subset B(b_j(h),h^\rho)$
for $j\in J_\natural\,$, and
$\chi_{j,h}\equiv 1$ (resp. $\zeta_{j,h}^\natural\equiv1$) on $\bar B(a_j(h),h^\rho/2)$ (resp. $\bar B(b_j(h),h^\rho/2)$)\,.\\
We set, for $j\in J_\natural(h)$, $\eta_{j,h}^\natural = \zeta_{j,h}^\natural\mathbf{1}_{\bar\Omega}\,$.\\
Notice that for all $\alpha\in\mathbb{N}^n\,$, 
\beq\label{supPaCutoff}
\sup|\pa^\alpha\chi_{j,h}| = \mathcal{O}(h^{-|\alpha|\rho})~~~\textrm{ and }~~~
\sup|\pa^\alpha\eta_{j,h}^\natural| =  \mathcal{O}(h^{-|\alpha|\rho})\,.
\eeq

Now we introduce our approximating operators. For $j\in J_{int}(h)\,$, we set
\beq\label{defCAjh}
\left\{\begin{array}{l}
\CA_{j,h} = -h^2\Delta+i(V(a_j(h))+\nabla V(a_j(h))\cdot(x-a_j(h)))\,,\\
\CD(\CA_{j,h}) = H^2(\mathbb{R}^n)\cap L^2(\mathbb{R}^n ; |x|^2dx)\,.
\end{array}\right.
\eeq
Then, according to Lemma \ref{lemModelRn} and by rescaling $x\mapsto h^{-2/3}x\,$,
we have $\sigma(\CA_{j,h}) = \emptyset\,$, and for all $\omega\in\mathbb{R}\,$, there exists $C_\omega^0>0$ such that
\beq\label{resCAjh}
\sup_{\Re z\leq\omega h^{2/3}}\|(\CA_{j,h}-z)^{-1}\|\leq \frac{C_\omega}{h^{2/3}}\,.
\eeq

In order to define the approximating operators at the boundary,
for $\natural = \perp$, $\sslash$, $\angle$ and $j\in J_\natural(h)$, we denote by
$\mathcal{F}_{b_j} = \mathcal{F}_{b_j(h)}$ the local diffeomorphism
defined by (\ref{defF}), where we choose $b = b_j(h)$ as base point, so that $\varphi(0) = b_j(h)\,$. In these coordinates, 
we define our local approximation for $\CA_h$ near $b_j(h)$ as
\beq\label{defTildeCAjh}
\tilde\CA_{j,h}^\natural = -h^2\Delta_{y,z}+i\left(V(b_j(h))+\sum_{i=1}^{n-1}J_i^{(j)}y_j + J_n^{(j)}z\right)\,,
\eeq
where, for all $j\in J_\natural(h)$ and $i=1,\dots,n-1\,$,
\beq\label{defAlphaBeta(j)}
J_i^{(j)} = J_i^{(j)}(h) = \nabla V(b_j(h))\cdot\pa_i\varphi(0)~~~~\textrm{ and }~~~~
J_n^{(j)} = J_n^{(j)}(h) = \nabla V(b_j(h))\cdot\vec\nu(0)\,.
\eeq
Notice that, if $j\in J_\sslash(h)\,$, then
$J_n^{(j)} = 0\,$, hence $\tilde\CA_{j,h}^\sslash$ has the same form as operator $\CA_\sslash^D$ studied in Subsection~\ref{pParallel}. 
Hence, according to Lemma \ref{lemModelParallel} and after rescaling $(y,z)\mapsto (h^{-2/3}y,h^{-2/3}z)\,$, for all $\omega\in\mathbb{R}\,$, there exists
$C_\omega^1>0$ such that
\beq\label{resTildeCAjhParallel}
\sup_{\Re z\leq\omega h^{2/3}}\|(\tilde\CA_{j,h}^\sslash-z)^{-1}\|\leq\frac{C_\omega^1}{h^{2/3}}\,.
\eeq
Similarly, if $j\in J_\perp(h)$, then $(J_1^{(j)},\dots,J_{n-1}^{(j)}) = 0$ and $\tilde\CA_{j,h}^\perp$ has the same form
as operator $\CA_\perp$ considered in Subsection \ref{pPerp}, with $$ J_n = J_n^{(j)} = |\nabla V(b_j(h))|\,.$$
 Hence, in view of Lemma \ref{lemModelPerp}, 
for all $\omega < |\mu_1||\nabla V(b_j(h))|^{2/3}/2\,$, there exists $C_\omega^2>0$ such that
\beq\label{resTildeCAjhPerp}
\sup_{\Re z\leq\omega h^{2/3}}\|(\tilde\CA_{j,h}^\perp-z)^{-1}\|\leq\frac{C_\omega^2}{h^{2/3}}\,.
\eeq
If $j\in J_\angle(h)\,$, then $\tilde\CA_{j,h}^\angle$ has the form (\ref{defAJtheta}), with $|J| = |\nabla V(b_j(h))|\,,$
\beq\label{defTheta}
 \cos\theta = \frac{J_n^{(j)}}{|\nabla V(b_j(h))|}~~~\textrm{ and }~~~
  \vec v = \frac{1}{\sin^2\theta}\left(\frac{J_1^{(j)}}{|\nabla V(b_j(h))|},\dots,\frac{J_{n-1}^{(j)}}{|\nabla V(b_j(h))|}\right)\,.
\eeq
Let us define, for $\delta>0\,$, the subset of the boundary
\[
 \pa\Omega_\perp^{(\delta)} = \{x\in\pa\Omega : d(x,\pa\Omega_\perp)\leq\delta\}\,.
\]
Then, for any fixed $\varepsilon>0$, there exists $\delta_0>0$ such that, for all $x\in\pa\Omega_\perp^{(\delta_0)}\,$,
\beq\label{conditionEps}
\frac{|\mu_1|}{2}J_m^{2/3}-\varepsilon \leq \frac{|\mu_1|}{2}|\nabla V(x)|^{2/3}-\frac{\varepsilon}{2}\,,
\eeq
where $J_m$ is defined by (\ref{defJm}).\\
On the other hand, there exists $\theta_0\in(0,\pi/2)$ such that, for all $x\in\pa\Omega\setminus\pa\Omega_\perp^{(\delta_0)}\,$, the angle
$\theta$ defined in (\ref{defTheta}) satisfies $\theta\in[\theta_0,\pi-\theta_0]$\,.\\
Thus, for all $j\in J_\angle(h)\,$, by using (\ref{resAngle1}) (with $\varepsilon/2$ instead of $\varepsilon$) and (\ref{conditionEps})
if $b_j(h)\in\pa\Omega_\perp^{(\delta_0)}\,$, or by using (\ref{resAngle2}) if $b_j(h)\in\pa\Omega\setminus\pa\Omega_\perp^{(\delta_0)}\,$, there exists
$C_\varepsilon^3>0$ such that, for all $h\in(0,h_0)\,$,
\beq\label{resTildeCAjhQcq}
\sup_{\Re z\leq(|\mu_1|J_m^{2/3}/2-\varepsilon)h^{2/3}}\|(\tilde\CA_{j,h}^\angle-z)^{-1}\|\leq\frac{C_\varepsilon^3}{h^{2/3}}\,,
\eeq
(here again, we used the rescaling $(y,z)\mapsto (h^{-2/3}y,h^{-2/3}z)$~)\,.\\

Now let us check that the potential in (\ref{defTildeCAjh}) is a good approximation of the potential $iV\circ\mathcal{F}_{b_j}(x)$ near $b_j(h)$.\\
As $(y,z)\rightarrow0$, we have
\[
\varphi(y,z) = b_j(h)+\sum_{i=1}^{n-1}\pa_i\varphi(0)y_i + \mathcal{O}(|y|^2)\,,
\]
hence
\[
 \mathcal{F}_{b_j}(y,z) = b_j(h)+\sum_{i=1}^{n-1}\pa_i\varphi(0)y_i +z\vec\nu(y) + \mathcal{O}(|y|^2)\,,
\]
and using (\ref{G0}), we get
\begin{eqnarray}
&& \nabla V(b_j(h))\cdot (\mathcal{F}_{b_j}(y,z)-b_j(h))) \nonumber \\
&& =\left(\sum_{i=1}^{n-1} J_i^{(j)}\pa_i\varphi(0) + J_n^{(j)}\vec\nu\right)
\cdot \left(\sum_{k=1}^{n-1}y_k\pa_k\varphi(0)+z\vec\nu(y) + \mathcal{O}(|y|^2)\right) \nonumber \\
&& = \sum_{i=1}^{n-1}J_i^{(j)}y_i+J_n^{(j)}z+\mathcal{O}(|y|^2)\,.
\end{eqnarray}
Thus, using that
\[
 V(x) = V(b_j(h))+\nabla V(b_j(h))\cdot(x-b_j(h))+\mathcal{O}(|x-b_j(h)|^2)\,,
\]
we obtain
\beq\label{approxV}
V\circ\mathcal{F}_{b_j} = V(b_j(h))+ \sum_{i=1}^{n-1}J_i^{(j)}y_i+J_n^{(j)}z+\mathcal{O}(|(y,z)|^2)\,.
\eeq
Since $\Omega$ is compact, there exists a fixed neighborhood of $\pa\Omega$ which is covered by a finite number of charts $(\mathcal{F},\omega)$ as
defined in (\ref{defF}). Hence, up to a translation there is a finite number of diffeomorphisms $\mathcal{F}_{b_j}$ for $h\in(0,h_0)$ and
$j\in J_{bdry}(h)$. Consequently all the remainder terms $\mathcal{O}(|y|^2)$ and $\mathcal{O}(|(y,z)|^2)$ above are uniform with respect
to $j$ and $h\in(0,h_0)$.
\\

Now we gather the resolvents of the approximate operators previously defined to build an approximate resolvent for $\CA_h$. For a fixed $\varepsilon>0$ and
any $\nu\in\mathbb{R}$, we set
\beq\label{defLambdah1}
\lambda(h) = \lambda_0h^{2/3}+i\nu\,,~~~\lambda_0 = \frac{|\mu_1|}{2}J_m^{2/3}-\varepsilon\,.
\eeq
Let $\psi_{j,h}\in\mathcal{C}_0^\infty(\omega\cap\bar\Omega)$ and $\tilde\psi_{j,h}\in\mathcal{C}_0^\infty(\mathcal{U}\times(0,z_0))$
such that $\psi_{j,h}(x) = 1$ near $b_j(h)$ and $\tilde\psi_{j,h}(y,z) = 1$ near $0$. Here $\omega$, $\mathcal{U}$ and $z_0$ are the objects
appearing in (\ref{transfoF}) corresponding to the diffeomorphism $\mathcal{F}_{b_j}$ near $b_j(h)$. Then we set
\beq\label{defResBord}
 R_{j,h}^\natural = T_{\mathcal{F}_{b_j}}^{-1}\tilde\psi_{j,h}(\tilde\CA_{j,h}^\natural-\lambda(h))^{-1}T_{\mathcal{F}_{b_j}}\psi_{j,h}\,,
\eeq
where $T_{\mathcal{F}_{b_j}} = T_{\mathcal{F}_{b_j(h)}}$ is defined in (\ref{transfoF}).\\
Now we define our global approximate resolvent, for $h\in(0,h_0)\,$, by
\begin{eqnarray}
 \mathcal{R}(h) & = & \sum_{j\in J_{int}(h)}\chi_{j,h}(\CA_{j,h}-\lambda(h))^{-1}\chi_{j,h} \nonumber \\ 
 & & + 
\sum_{\tiny{\natural \in\{\perp,\sslash,\angle\}}~}
 \sum_{j\in J_\natural(h)}\eta_{j,h}^\natural R_{j,h}^\natural\eta_{j,h}^\natural\,.
 \label{defResApp1}
\end{eqnarray}
Then, we have
\begin{eqnarray}
 (\CA_h-\lambda(h))\mathcal{R}(h) &= &I  +  \sum_{j\in J_{int}(h)}\chi_{j,h}(\CA_h-\CA_{j,h})(\CA_{j,h}-\lambda(h))^{-1}\chi_{j,h} \nonumber \\
 & + & \sum_{j\in J_{int}(h)}[\CA_h,\chi_{j,h}](\CA_{j,h}-\lambda(h))^{-1}\chi_{j,h} \nonumber \\
 & + & \sum_{\tiny{\natural \in\{\perp,\sslash,\angle\}}~}\sum_{j\in J_\natural(h)}
 T_{{\mathcal{F}_{b_j}}}^{-1}\tilde\eta_{j,h}^\natural(\tilde\CA_h-\tilde\CA_{j,h}^\natural)(\tilde\CA_{j,h}^\natural-\lambda(h))^{-1}
 T_{\mathcal{F}_{b_j}}\eta_{j,h}^\natural\nonumber\\
  & + & \sum_{\tiny{\natural \in\{\perp,\sslash,\angle\}}~}\sum_{j\in J_\natural(h)}
  T_{\mathcal{F}_{b_j}}^{-1}[\tilde\CA_h,\tilde\eta_{j,h}^\natural](\tilde\CA_{j,h}^\natural-\lambda(h))^{-1}T_{\mathcal{F}_{b_j}}\eta_{j,h}^\natural\,,
  \label{decompRlambda}
\end{eqnarray}
where $\tilde\CA_h = T_{\mathcal{F}_{b_j}}\psi_{j,h}\CA_hT_{\mathcal{F}_{b_j}}^{-1}\tilde\psi_{j,h}$ denotes the operator $\CA_h$ expressed in the local coordinates
near $b_j(h)$ (see (\ref{AhCoordLocales})), and $\tilde\eta_{j,h}^\natural = \eta_{j,h}^\natural\circ{\mathcal{F}_{b_j}}$.\\
In the following, we estimate each term of the right-hand side.\\

First, for $j\in J_{int}(h)$, we have
\[
 \CA_h-\CA_{j,h} = i\mathcal{O}(|x-a_j(h)|^2),~~~x\rightarrow a_j(h),
\]
hence $\|\chi_{j,h}(\CA_h-\CA_{j,h})\| = \mathcal{O}(h^{2\rho})$. According to (\ref{resCAjh}), we then get
\beq\label{estRHS_R1}
\|\chi_{j,h}(\CA_h-\CA_{j,h})(\CA_{j,h}-\lambda(h))^{-1}\chi_{j,h}\| = \mathcal{O}(h^{2(\rho-1/3)})\,.
\eeq
Now we estimate the terms of the second sum in the right-hand side of (\ref{decompRlambda}). We have, for $j\in J_{int}(h)\,$,
\begin{eqnarray}
 [\CA_h,\chi_{j,h}](\CA_{j,h}-\lambda(h))^{-1}\chi_{j,h} &= & 
 -h^2\Delta\chi_{j,h}(\CA_{j,h}-\lambda(h))^{-1}\chi_{j,h} \nonumber \\
  & & - 2h\nabla\chi_{j,h}\cdot h\nabla(\CA_{j,h}-\lambda(h))^{-1}\chi_{j,h} \nonumber \\
   & =: & P_{j,h}^{(1)}+P_{j,h}^{(2)}\,.\label{decompCommut1}
 \end{eqnarray}
According to (\ref{supPaCutoff}) and (\ref{resCAjh}),
\beq\label{estP1}
\|P_{j,h}^{(1)}\| = \mathcal{O}(h^{2(2/3-\rho)})\,.
\eeq
On the other hand, for every $v\in\CD(\CA_{j,h})$, we have
\begin{eqnarray*}
\|h\nabla v\|^2 & = & \Re\sca{(\CA_{j,h}-\lambda(h))v}{v} + \Re\lambda(h)\|v\|^2 \\
 & \leq & \|(\CA_{j,h}-\lambda(h))v\|\|v\|+\Re\lambda(h)\|v\|^2\,,
 \end{eqnarray*}
which, applied to $v = (\CA_{j,h}-\lambda(h))^{-1}\chi_{j,h}f$, $f\in L^2(\mathbb{R}^n)$, yields
\begin{eqnarray}
 \|h\nabla(\CA_{j,h}-\lambda(h))^{-1}\chi_{j,h}f\|  \leq \sqrt{\|(\CA_{j,h}-\lambda(h))^{-1}\|}\|f\| &&\nonumber \\
  +\sqrt{|\Re\lambda(h)|}\|(\CA_{j,h}-\lambda(h))^{-1}\|\|f\|\,,&& \label{estApriori}
\end{eqnarray}
that is, in view of (\ref{supPaCutoff}), (\ref{resCAjh}) and (\ref{defLambdah1})\,,
\beq\label{estP2}
 \|P_{j,h}^{(2)}\| = \mathcal{O}(h^{2/3-\rho})\,.
 \eeq
Thus, (\ref{decompCommut1}), (\ref{estP1}) and (\ref{estP2}) yield, for every $j\in J_{int}(h)\,$,
\beq\label{estCommut1}
\|[\CA_h,\chi_{j,h}](\CA_{j,h}-\lambda(h))^{-1}\chi_{j,h}\| = \mathcal{O}(h^{2/3-\rho})\,.
\eeq
\\
Now we consider the boundary terms in (\ref{decompRlambda}). First, according to (\ref{AhCoordLocales}) and (\ref{approxV}), for 
$\natural = \sslash, \perp, \angle$ and $j\in J_\natural(h)\,$, we have
\[
\tilde\CA_h-\tilde\CA_{j,h}^\natural = -h^2\sum_{i,k=1}^{n-1}(G^{ik}(y,z)-1)\pa_i\pa_k-h^2\sum_{i=1}^n\beta_i(y,z)\pa_i + \mathcal{O}(|(y,z)|^2)\,.
\]
Here the functions $G^{ik}$ and $\beta_i$ depend on the index $j\in J_\natural(h)\,$, although we do not mention it in the notation. However,
the remainder term $\mathcal{O}(|(y,z)|^2)$ is uniform with respect to $j$ and $h\in(0,h_0)$.\\
Hence, according
to (\ref{DLGij}), and since the functions $\beta_i$ are bounded on $\Supp\tilde\eta_{j,h}^\natural$, we have, for some $C>0\,$,
\begin{eqnarray}
& & \|\tilde\eta_{j,h}^\natural(\tilde\CA_h-\tilde\CA_{j,h}^\natural)(\tilde\CA_{j,h}^\natural-\lambda(h))^{-1}\|
\leq Ch^\rho\|h^2\Delta_{(y,z)}(\tilde\CA_{j,h}^\natural-\lambda(h))^{-1}\| \nonumber \\
 & & + C\Big(\|h^2\nabla_{(y,z)}(\tilde\CA_{j,h}^\natural-\lambda(h))^{-1}\| +
 h^{2\rho}\|(\tilde\CA_{j,h}^\natural-\lambda(h))^{-1}\|\Big)\,.\label{decompDifferenceBdry}
\end{eqnarray}
Regarding the second term in the right-hand side of (\ref{decompDifferenceBdry}), we can use an estimate similar to (\ref{estApriori})
to get
\[
\|h\nabla_{(y,z)}(\tilde\CA_{j,h}^\natural-\lambda(h))^{-1}\|\leq \|(\tilde\CA_{j,h}^\natural-\lambda(h))^{-1}\|^{1/2} +
(\Re\lambda(h))^{1/2}\|(\tilde\CA_{j,h}^\natural-\lambda(h))^{-1}\|\,,
\]
hence using (\ref{resTildeCAjhParallel}), (\ref{resTildeCAjhPerp}) and (\ref{resTildeCAjhQcq}),
\beq\label{DiffBdry2}
\|h\nabla_{(y,z)}(\tilde\CA_{j,h}^\natural-\lambda(h))^{-1}\| = \mathcal{O}(h^{-1/3})\,.
\eeq

The first norm in the right-hand side can be estimated as follows: as in (\ref{estDomA+}), we can write,
for all $u\in\CD(\tilde\CA_{j,h}^\natural)\,$,
\begin{eqnarray}
 \|(\tilde\CA_{j,h}^\natural-\lambda(h))u\|^2 & = & \|h^2\Delta u\|^2 + 
 \|(i\tilde\ell-\lambda(h))u\|^2  + 2\Re\sca{-h^2\Delta u}{(i\tilde\ell-\lambda(h))u} \nonumber \\
&\geq &\|h^2\Delta u\|^2 + 2\Re\sca{h\nabla u}{h\nabla[(i\tilde\ell-\lambda(h))u]} \nonumber \\
 & \geq &\|h^2\Delta u\|^2+ 2h\Im\sca{h\nabla u}{J^{(j)}u} -2\lambda_0h^{2/3}\|h\nabla u\|^2\,,
 \label{DiffBdrInter}
\end{eqnarray}
where $J^{(j)} = (J_1^{(j)},\dots,J_n^{(j)})\,$, and $\tilde\ell$ denotes the potential of $\tilde\CA_{j,h}^\natural\,$,
\[
V(b_j(h))+J_1^{(j)}y_1+\dots+ J_n^{(j)}z\,.
\]
Notice that $|J^{(j)}|$ is bounded uniformly with respect to $j$ and $h$ since $\pa\Omega$ is compact.\\
Thus, for some $c>0$ independent of $j$ and $h\,$,
\beq
\|h^2\Delta u\|^2 \leq \|(\tilde\CA_{j,h}^\natural-\lambda(h))u\|^2
+ c\,(h\|h\nabla u\|\|u\| + h^{2/3}\|h\nabla u\|^2)\,.
\eeq 
Applying this estimate to $u = (\tilde\CA_{j,h}^\natural-\lambda(h))^{-1}f$, with $f\in L^2(\mathbb{R}_+^n)$, we then get
from (\ref{resTildeCAjhParallel}), (\ref{resTildeCAjhPerp}), (\ref{resTildeCAjhQcq}) and (\ref{DiffBdry2}),
\[
 \|h^2\Delta(\tilde\CA_{j,h}^\natural-\lambda(h))^{-1}f\|^2 \leq c'\,\|f\|^2\,,~~~~c'>0\,,
\]
that is
\beq\label{DiffBdry1}
h^\rho\|\tilde\eta_{j,h}^\natural h^2\Delta_{(y,z)}(\tilde\CA_{j,h}^\natural-\lambda(h))^{-1}\| = \mathcal{O}(h^\rho)\,.
\eeq
Then, (\ref{decompDifferenceBdry}), (\ref{DiffBdry2}) and (\ref{DiffBdry1}) yield
\beq\label{estRHS_R1bdry}
\|T_{{\mathcal{F}_{b_j}}}^{-1}\tilde\eta_{j,h}^\natural(\tilde\CA_h-\tilde\CA_{j,h}^\natural)(\tilde\CA_{j,h}^\natural-\lambda(h))^{-1}
 T_{\mathcal{F}_{b_j}}\eta_{j,h}^\natural\| = \mathcal{O}(h^\rho) + \mathcal{O}(h^{2/3}) + \mathcal{O}(h^{2(\rho-1/3)})\,.
\eeq
\\
Finally, the terms contained in the last sum of the right-hand side in (\ref{decompRlambda}) can be estimated as in (\ref{estCommut1}): 
\beq\label{estCommut1bdry}
\|T_{\mathcal{F}_{b_j}}^{-1}[\tilde\CA_h,\tilde\eta_{j,h}^\natural](\tilde\CA_{j,h}^\natural-\lambda(h))^{-1}T_{\mathcal{F}_{b_j}}\eta_{j,h}^\natural\| 
= \mathcal{O}(h^{2/3-\rho})\,.
\eeq
Let us stress that the constants in estimates (\ref{resCAjh}), (\ref{resTildeCAjhParallel}), (\ref{resTildeCAjhPerp}) and (\ref{resTildeCAjhQcq})
are independent of $h\in(0,h_0)$ and $j=j(h)$. Hence estimates (\ref{estRHS_R1}), (\ref{estCommut1}), (\ref{estRHS_R1bdry}) and
(\ref{estCommut1bdry}) are uniform with respect to $j$.\\

For the time being, we have controlled separately each term appearing in the right-hand side of (\ref{decompRlambda}). However, since
the sums therein contain a growing number of terms as $h \rightarrow 0\,$, we shall sum these estimates carefully 
in order to get an appropriate
bound eventually. In this purpose, we take into account the {\em almost orthogonality} of those terms. Namely, we use the following lemma
(\cite{Stein}, VII \S $2$):
\begin{lemma}[Cotlar-Stein Lemma]\label{lemCotlar}
Let $\CH$ be a Hilbert space and, for every $j\in\mathbb{N}\,$, $T_j\in\mathcal{L}(\CH)\,$. Assume that
 \beq\label{hypCotlarA}
 A :=\sup_{j\in\mathbb{N}}\sum_{k\in\mathbb{N}}\sqrt{\|T_jT_k^*\|} < +\infty\,,
 \eeq
 and
  \beq\label{hypCotlarB}
 B :=\sup_{j\in\mathbb{N}}\sum_{k\in\mathbb{N}}\sqrt{\|T_j^*T_k\|} < +\infty\,.
 \eeq
Then, $\sum_{j\in\mathbb{N}}T_j$ converges in the strong operator topology and
 \beq\label{cclCotlar}
 \left\|\sum_{j\in\mathbb{N}}T_j\right\|\leq\sqrt{AB}\,.
 \eeq
\end{lemma}

Here we apply this lemma for a fixed $h\in(0,h_0)$ with, for $j\in J_{int}(h)\,$, 
\[
 T_j = \chi_{j,h}(\CA_h-\CA_{j,h})(\CA_{j,h}-\lambda(h))^{-1}\chi_{j,h}\,.
\]
Since the first sum in the right-hand side of (\ref{decompRlambda}) has a finite number of terms, we set
$T_j = 0$ for $j\in\mathbb{N}\setminus J_{int}(h)$.
By definition of the functions $\chi_{j,h}$, for every $j_0\in J_{int}(h)\,$, 
we have $\Supp\chi_{j_0,h}\cap\Supp\chi_{k,h}=\emptyset$ except for a finite number (uniformly bounded with respect to $h$ and $j_0$) of indices $k\,$,
which we shall denote by $\{k_1(j_0),\dots,k_p(j_0)\}\,$.\\
 Hence we have clearly
\[
 \forall j_0\in J_{int}(h)\,,~k\in J_{int}(h)\setminus\{k_1(j_0),\dots,k_p(j_0)\}\,,~~~T_{j_0}T_k^* = T_{j_0}^*T_k = 0\,.
\]
Thus,
\[
 A = \sup_{j_0\in J_{int}(h)}\sum_{\ell=1}^p\|T_{j_0}T_{k_\ell(j_0)}^*\| = \mathcal{O}(h^{2(\rho-1/3)})\,,
\]
according to (\ref{estRHS_R1}).\\
Similarly, we get $B= \mathcal{O}(h^{2(\rho-1/3)})$\,. Hence (\ref{cclCotlar}) yields
\beq\label{estSommeR1}
\left\|\sum_{j\in J_{int}(h)}\chi_{j,h}(\CA_h-\CA_{j,h})(\CA_{j,h}-\lambda(h))^{-1}\chi_{j,h}\right\| = \mathcal{O}(h^{2(\rho-1/3)})\,.
\eeq
We can handle the other sums in (\ref{decompRlambda}) alike to get, in view of (\ref{estCommut1}), (\ref{estRHS_R1bdry}) and (\ref{estCommut1bdry}),
\beq\label{estSommeR2}
\left\|\sum_{j\in J_{int}(h)}[\CA_h,\chi_{j,h}](\CA_{j,h}-\lambda(h))^{-1}\chi_{j,h}\right\| = \mathcal{O}(h^{2/3-\rho})\,,
\eeq
\beq\label{estSommeR3}
\left\|\sum_{\tiny{\natural \in\{\perp,\sslash,\angle\}}~}\sum_{j\in J_\natural(h)}
 T_{{\mathcal{F}_{b_j}}}^{-1}\tilde\eta_{j,h}^\natural(\tilde\CA_h-\tilde\CA_{j,h}^\natural)(\tilde\CA_{j,h}^\natural-\lambda(h))^{-1}
 T_{\mathcal{F}_{b_j}}\eta_{j,h}^\natural\right\| = \mathcal{O}(h^\rho) + \mathcal{O}(h^{2/3}) + \mathcal{O}(h^{2(\rho-1/3)})\,,
 \eeq
 and
 \beq\label{estSommeR4}
 \left\|\sum_{\tiny{\natural \in\{\perp,\sslash,\angle\}}~}\sum_{j\in J_\natural(h)}
  T_{\mathcal{F}_{b_j}}^{-1}[\tilde\CA_h,\tilde\eta_{j,h}^\natural](\tilde\CA_{j,h}^\natural-\lambda(h))^{-1}T_{\mathcal{F}_{b_j}}\eta_{j,h}^\natural
  \right\| = \mathcal{O}(h^{2/3-\rho})\,.
\eeq
Thus, if we choose $\rho\in(1/3,2/3)\,$, we obtain from (\ref{decompRlambda}), (\ref{estSommeR1}), (\ref{estSommeR2}), (\ref{estSommeR3}) and
(\ref{estSommeR4}):
\beq\label{estReste1}
(\CA_h-\lambda(h))\CR(h) = I+\CE(h)\,,~~~~\|\CE(h)\|\underset{\tiny{h\rightarrow0}}{\longrightarrow}0\,.
\eeq
Hence, there exists $h_\varepsilon\in(0,h_0)$ such that, for all $h\in(0,h_\varepsilon)\,$, $(\CA_h-\lambda(h))$ is invertible, with
\[
 (\CA_h-\lambda(h))^{-1} = \CR(h)(I+\CE(h))^{-1}\,.
\]
Consequently, there is a strip free from eigenvalues:
\[
 \forall h\in(0,h_\varepsilon)\,,~~~\sigma(\CA_h)\cap([0,(|\mu_1|J_m^{2/3}/2-\varepsilon)h^{2/3}]+
 i\mathbb{R}) = \emptyset\,,
\]
which proves (\ref{limSpect1}).
Moreover, we have of course $\|(I+\CE(h))^{-1}\| = \mathcal{O}(1)\,$,
and according to (\ref{resCAjh}), (\ref{resTildeCAjhParallel}), (\ref{resTildeCAjhPerp}) and (\ref{resTildeCAjhQcq}), by using Lemma \ref{lemCotlar} again to estimate the
sums in (\ref{defResApp1}), we get
\[\|\CR(h)\| = \mathcal{O}(h^{-2/3})\,.\]
The estimate (\ref{estRes1}) follows, which concludes the proof of Theorem \ref{thmNSAschro1}.

\section{Lower bound for a Morse potential}\label{sthm3Lower}
Here we prove part of the statements in Theorem \ref{thmNSAschro2}. Namely, we prove (\ref{estRes2}), as well as the lower bound in
(\ref{limSpect2}):
\beq\label{LowerlimSpect2}
\varliminf_{h\rightarrow0}\frac{1}{h}\inf\Re\sigma(\CA_h) \geq \frac{\kappa}{2}\,.
\eeq
The corresponding upper bound shall be proved in Section \ref{sthm3Upper}.\\

We follow the same method as in Section \ref{sthm1}, but we will need a quadratic approximation in the neighborhood of the critical points of $V$.\\
Let $x_1^c,\dots,x_p^c$ be the critical points of $V\,$, and, for $k=1,\dots,p\,$, $\theta_{k,h}\in\mathcal{C}_0^\infty(\mathbb{R}^n ; [0,1])\,$
such that 
\beq\label{SupportTheta}
\Supp\theta_{k,h}\subset\{x\in\mathbb{R}^n : |x-x_k^c|\leq h^{\rho'}\}\,,
\eeq
where $\rho'>0$ shall be determined later, and
$\theta_{k,h}(x) = 1$ if $|x-x_k^c|\leq h^{\rho'}/2\,$. As in Section \ref{sthm1}, for any $h\in(0,h_0)$ we consider a covering
of the compact set $\bar\Omega\setminus\bigcup_{k=1}^p\{|x-x_k^c|> h^{\rho'}\}$ by balls 
$B(a_j(h),h^\rho/2)$ and $B(b_k(h),h^\rho/2)\,$, $j\in J_{int}(h)\,$, $k\in J_{bdry}(h)\,$, where $a_j(h)\in\Omega$ and
$b_j(h)\in\pa\Omega\,$, such that the corresponding closed balls
of radius $h^\rho/2$ do not intersect one another, and such that, for every $h\in(0,h_0)\,$, $j\in J_{int}(h)$ and $k=1,\dots,p\,$,
\beq\label{BoulesIntersVide}
\bar B(a_j(h),h^\rho/2)\cap\bar B(x_k^c,h^{\rho'}/2) = \emptyset\,.
\eeq
Then we define $J_\natural(h)\,$, $\natural = \sslash,
\perp,\angle\,$, as in Section \ref{sthm1}, as well as the functions $\chi_{j,h}$ and $\eta_{j,h}^\perp\,$, with the following condition instead of
(\ref{sommePartition1}):
\beq\label{sommePartition2}
\forall x\in\bar\Omega\,,~~~
\sum_{j\in J_{int}(h)}\chi_{j,h}(x)^2+\sum_{\natural,~k\in J_\natural(h)}\eta_{k,h}^\natural(x)^2+\sum_{k=1}^p\theta_{k,h}(x)^2=1\,.
\eeq
Let $\CA_{j,h}\,$, $\tilde\CA_{j,h}^\natural(h)$ and $R_{j,h}^\natural$ denote the same approximate operators as before. For $k=1,\dots,p\,$, we set
\beq\label{defHkh}
\left\{\begin{array}{l}
\CH_{k,h} = -h^2\Delta+i(V(x_k^c)+(\Hess V(x_k^c)(x-x_k^c))\cdot(x-x_k^c))\,, \\
\CD(\CH_{k,h}) = H^2(\mathbb{R}^n)\cap L^2(\mathbb{R}^n ; |x|^4dx)\,,
\end{array}\right.
\eeq
which will stand as an approximation of $\CA_h$ near $x_k^c$\,.\\
Instead of (\ref{defLambdah1}) we set, for any $\nu\in\mathbb{R}$ and $\varepsilon>0\,$,
\beq\label{defLambdah2}
\lambda(h) = \lambda_0h+i\nu\,,~~~\lambda_0 = \frac{\kappa}{2}-\varepsilon\,,
\eeq
where $\kappa$ is the constant in Theorem \ref{thmNSAschro2}.\\
Our approximate resolvent will be
\begin{eqnarray}
 \mathcal{Q}(h) & = & \sum_{j\in J_{int}(h)}\chi_{j,h}(\CA_{j,h}-\lambda(h))^{-1}\chi_{j,h} \nonumber \\ 
 & & + \sum_{\tiny{\natural \in\{\perp,\sslash,\angle\}}~}
 \sum_{j\in J_\natural(h)}\eta_{j,h}^\natural R_{j,h}^\natural\eta_{j,h}^\natural \nonumber \\
  & & +\sum_{k=1}^p\theta_{k,h}(\CH_{k,h}-\lambda(h))^{-1}\theta_{k,h}\,.
 \label{defResApp2}
\end{eqnarray}
Then, we have
\begin{eqnarray}
 (\CA_h-\lambda(h))\mathcal{Q}(h) &=& I + \sum_{j\in J_{int}(h)}\chi_{j,h}(\CA_h-\CA_{j,h})(\CA_{j,h}-\lambda(h))^{-1}\chi_{j,h} \nonumber \\
 && +  \sum_{j\in J_{int}(h)}[\CA_h,\chi_{j,h}](\CA_{j,h}-\lambda(h))^{-1}\chi_{j,h} \nonumber \\
 && +  \sum_{\tiny{\natural \in\{\perp,\sslash,\angle\}}~}\sum_{j\in J_\natural(h)}
 T_{{\mathcal{F}_{b_j}}}^{-1}\tilde\eta_{j,h}^\natural(\tilde\CA_h-\tilde\CA_{j,h}^\natural)(\tilde\CA_{j,h}^\natural-\lambda(h))^{-1}
 T_{\mathcal{F}_{b_j}}\eta_{j,h}^\natural \nonumber\\
  && +  \sum_{\tiny{\natural \in\{\perp,\sslash,\angle\}}~}\sum_{j\in J_\natural(h)}
  T_{\mathcal{F}_{b_j}}^{-1}[\tilde\CA_h,\tilde\eta_{j,h}^\natural](\tilde\CA_{j,h}^\natural-\lambda(h))^{-1}T_{\mathcal{F}_{b_j}}\eta_{j,h}^\natural \nonumber \\
  && +  \sum_{k=1}^p\theta_{k,h}(\CA_h-\CH_{k,h})(\CH_{k,h}-\lambda(h))^{-1}\theta_{k,h} \nonumber \\
   && + \sum_{k=1}^p[\CA_h,\theta_{k,h}](\CH_{k,h}-\lambda(h))^{-1}\theta_{k,h}\,,
  \label{decompQlambda}
\end{eqnarray}
where $\tilde\CA_h = T_{\mathcal{F}_{b_j}}\psi_{j,h}\CA_hT_{\mathcal{F}_{b_j}}^{-1}\tilde\psi_{j,h}$ denotes the operator $\CA_h$ expressed in the local coordinates
near $b_j(h)$ (see (\ref{AhCoordLocales})), and $\tilde\eta_{j,h}^\natural = \eta_{j,h}^\natural\circ{\mathcal{F}_{b_j}}$\,.\\

The boundary terms, that is those appearing in the third and fourth sums in the righ-hand side, can be estimated as in Section \ref{sthm1}, hence
(\ref{estSommeR3}) and (\ref{estSommeR4}) hold.\\

Regarding the first and second sums in the right-hand side, we have to take into account that, when $a_j(h)$ is close to some critical point $x_k^c\,$,
$|\nabla V(a_j(h))|$ can become small as $h\rightarrow0$. However, according to (\ref{BoulesIntersVide}), we have for all $j\in J_{int}(h)\,$,
\[
\forall x\in\Supp\chi_{j,h}\,,~~~|x-x_k^c|\geq\frac{h^{\rho'}}{2}\,.
\]
Hence, using that
\[
 \nabla V(a_j(h)) = \Hess V(x_k^c)\cdot(a_j(h)-x_k^c) + \mathcal{O}(|a_j(h)-x_k^c|^2)\,,
\]
there exists $c>0$ such that, for every $h\in(0,h_0)$ and $j\in J_{int}(h)$\,,
\beq\label{borneInfNabla}
|\nabla V(a_j(h))|\geq \frac{h^{\rho'}}{c}\,.
\eeq
According to Subsection \ref{pWhole}, after a rotation we can assume that $\CA_{j,h}$ has the form
\[
-h^2\Delta+i|\nabla V(a_j(h)|x_1 + i(V(a_j(h))-\nabla V(a_j(h))\cdot a_j(h))\,.
\]
Now if $T_h$ denotes the unitary map 
\[T_h : u(x)\mapsto \frac{h^{2/3}}{|\nabla V(a_j(h))|^{1/3}}u\left(\frac{|\nabla V(a_j(h))|^{1/3}}{h^{2/3}}x\right)\,,\]
then
\[
\left\{
\begin{array}{l}
 \CA_{j,h}-\lambda(h) = \big(h|\nabla V(a_j(h))|\big)^{2/3}T_h^{-1}\Big(\CA_0-\big(h|\nabla V(a_j(h))|\big)^{-2/3}\big(\lambda(h)+i\nu_0(h)\big)
 \Big)T_h\,, \\
 \nu_0(h) = V(a_j(h))-\nabla V(a_j(h))\cdot a_j(h)\,,
 \end{array}
 \right.
\]
where $\CA_0 = -\Delta+ix_1$ is the operator of Subsection \ref{pWhole}. \\
Thus, we have
\beq\label{resAjhinter}
 \|(\CA_{j,h}-\lambda(h))^{-1}\| = \frac{1}{\big(h|\nabla V(a_j(h))|\big)^{2/3}}
 \Big\|\Big(\CA_0-\big(h|\nabla V(a_j(h))|\big)^{-2/3}\big(\lambda(h)+i\nu_0(h)\big)\Big)^{-1}\Big\|\,.
\eeq
Besides, in view of (\ref{defLambdah2}) and (\ref{borneInfNabla}), if we choose $\rho'<1/2\,$, then there exists $\omega>0$ such that, for all
$h\in(0,h_0)\,$,
\[
 (h|\nabla V(a_j(h))|)^{-2/3}\Re\lambda(h) \leq\omega\,.
\]
Hence, (\ref{borneInfNabla}), (\ref{resAjhinter}) and Lemma \ref{lemModelRn} yield
\beq\label{resAjh}
 \|(\CA_{j,h}-\lambda(h))^{-1}\| = \mathcal{O}\left(\frac{1}{h^{2/3(1+\rho')}}\right)\,.
 \eeq
Using this resolvent estimate, we prove as for (\ref{estRHS_R1}) that, for all $j\in J_{int}(h)\,$,
\beq\label{estRHS_Q1}
\|\chi_{j,h}(\CA_h-\CA_{j,h})(\CA_{j,h}-\lambda(h))^{-1}\chi_{j,h} \| = \mathcal{O}(h^{2\rho-2/3(1+\rho')})\,.
\eeq
\\

Now we handle the commutator terms as in Section \ref{sthm1}, by estimating the two terms of
(\ref{decompCommut1}). First (\ref{estP1}) clearly becomes
\beq\label{estP1Q}
\|P_{j,h}^{(1)}\| = \mathcal{O}(h^{4/3-2\rho-2\rho'/3})\,.
\eeq
On the other hand, (\ref{estApriori}) and (\ref{resAjh}) imply
\begin{eqnarray*}
 \|h\nabla(\CA_{j,h}-\lambda(h))^{-1}\chi_{j,k}\|& = &\mathcal{O}\left(\frac{1}{h^{1/3(1+\rho')}}\right) +
 \mathcal{O}\left(\frac{1}{h^{1/6+2\rho'/3}}\right)\\
  & = & \mathcal{O}\left(\frac{1}{h^{1/3(1+\rho')}}\right)\,,
\end{eqnarray*}
since $\rho'<1/2\,$.\\
 Hence
\beq\label{estP2Q}
\|P_{j,h}^{(2)}\| = \mathcal{O}(h^{2/3-\rho-\rho'/3})\,,
\eeq
and by (\ref{estP1Q}) we get, for all $j\in J_{int}(h)\,$,
\beq\label{estCommut2}
\|[\CA_h,\chi_{j,h}](\CA_{j,h}-\lambda(h))^{-1}\chi_{j,h}\| = \mathcal{O}(h^{2/3-\rho-\rho'/3})\,.
\eeq
\\

It remains to estimate the terms of the two last sums in the right-hand side of (\ref{decompQlambda}). For each $k=1,\dots,p\,$, let $U_k$
be an orthogonal matrix such that
\[
 {}^tU_k~\Hess V(x_k^c)~U_k = \left(\begin{array}{ccc}
                                 \lambda_1^k & 0 & 0 \\
                                 0 & \ddots & 0 \\
                                 0 & 0 & \lambda_n^k
                                \end{array}
\right)\,,
\]
where $\{\lambda_j^k\}_{j=1,\dots,n} = \sigma(\Hess V(x_k^c))\,$.\\
Let $T_{h,k} : u(x)\mapsto u(h^{-1/2}U_k(x-x_k^c))\,$. Then,
\beq\label{rescaleQuad}
 T_{h,k}(\CH_{k,h}-\lambda(h))T_{h,k}^{-1} = h\left(-\Delta+\frac{i}{2}\sum_{j=1}^n\lambda_j^kx_j^2-(\lambda(h)-iV(x_k^c))h^{-1}\right)\,.
\eeq
Since $x_k^c$, $k=1\dots,n\,$, are non-degenerate critical points, we have $\lambda_j^k\neq0$ for $j=1,\dots,n\,$. Hence, according to Lemma
\ref{lemModelQuad}, we have
\[
 \inf\Re\sigma(\CH_{k,h}) = \frac{\kappa_k}{2}h\,,
\]
where $\kappa_k$ is the constant defined in (\ref{defKappak}). Moreover, since $\Re (\lambda(h)-iV(x_k^c))h^{-1}<\kappa_k/2$ for any $k=1,\dots,p\,$
due to (\ref{defLambdah2}), (\ref{rescaleQuad}) and (\ref{estResModelQuad}) yield
\beq\label{estResQuad}
\|(\CH_{k,h}-\lambda(h))^{-1}\| = \mathcal{O}\left(\frac{1}{h}\right)\,.
\eeq
On the other hand, according to (\ref{SupportTheta}),
\beq\label{borneSupTheta}
\forall k=1,\dots,p\,,~\forall \alpha\in\mathbb{N},~~~\sup|\pa^\alpha\theta_{k,h}| = \mathcal{O}(h^{-|\alpha|\rho'})\,,
\eeq
and
\beq\label{difference}
\|\theta_{k,h}(\CA_h-\CH_{k,h})\| = \mathcal{O}(h^{3\rho'})\,.
\eeq
Thus, combining (\ref{estResQuad}), (\ref{borneSupTheta}), (\ref{difference}), and following the proof of (\ref{estRHS_R1}) and (\ref{estCommut1}),
we get, for all $k=1,\dots,p\,$,
\beq\label{estRHS_Qquad}
\|\theta_{k,h}(\CA_h-\CH_{k,h})(\CH_{k,h}-\lambda(h))^{-1}\theta_{k,h}\| = \mathcal{O}(h^{3\rho'-1})
\eeq
and
\beq\label{estCommut2quad}
\|[\CA_h,\theta_{k,h}](\CH_{k,h}-\lambda(h))^{-1}\theta_{k,h}\| = \mathcal{O}(h^{1/2-\rho'})\,.
\eeq
\\

As a conclusion, if we choose 
\beq\label{rhorho'}  \frac{1}{3}<\rho'<\frac{1}{2}~~~~~\textrm{ and }~~~~~
  \frac{1+\rho'}{3}<\rho<\frac{2-\rho'}{3}\,,
\eeq
then according to (\ref{estSommeR3}), (\ref{estSommeR4}), (\ref{decompQlambda}), (\ref{estRHS_Q1}), (\ref{estCommut2}), (\ref{estRHS_Qquad}),
(\ref{estCommut2quad}), and by using again Lemma \ref{lemCotlar} to handle the large number of terms in the sums, we get
\beq\label{estReste2}
(\CA_h-\lambda(h))\mathcal{Q}(h) = I+\tilde\CE(h)\,,~~~~\|\tilde\CE(h)\|\underset{\tiny{h\rightarrow0}}{\longrightarrow}0\,.
\eeq
Hence, there exists $h_\varepsilon\in(0,h_0)$ such that, for all $h\in(0,h_\varepsilon)\,$, $(\CA_h-\lambda(h))$ is invertible, with
\[
 (\CA_h-\lambda(h))^{-1} = \mathcal{Q}(h)(I+\tilde\CE(h))^{-1}\,.
\]
Consequently, there is a strip free from eigenvalues:
\[
 \forall h\in(0,h_\varepsilon),~~~\sigma(\CA_h)\cap([0,(\kappa/2-\varepsilon)h]+i\mathbb{R}) = \emptyset\,,
\]
which proves (\ref{LowerlimSpect2}). Moreover, we have of course $\|(I+\CE(h))^{-1}\| = \mathcal{O}(1)\,$,
and according to (\ref{resTildeCAjhParallel}), (\ref{resTildeCAjhPerp}), (\ref{resTildeCAjhQcq}), (\ref{resAjh}), (\ref{estResQuad}), and by using Lemma \ref{lemCotlar} again to estimate the
sums in (\ref{defResApp2}), we get
\[\|\mathcal{Q}(h)\| = \mathcal{O}(h^{-1})\,.\]
The estimate (\ref{estRes1}) follows.

\section{Upper bound for a potential without critical point in dimension $1$}\label{sthm2}
\sectionmark{Dimension $1$ : upper bound}
In this section, we prove Theorem \ref{thmNSAschro1_1d}. In view of the statement of Theorem \ref{thmNSAschro1}, it only remains to prove that
\beq\label{upper1d}
\varlimsup\limits_{h\to0}\frac{1}{h^{2/3}}\inf\Re\sigma(\CA_h) \leq \frac{|\mu_1|}{2}J^{2/3}\,.
\eeq

Up to a scale change, we can assume that $a=0$ and $b=1\,$. Moreover, without loss of generality, we shall assume in this section that
$V'>0$ on $(0,1)$, and $J = |V'(0)|^{2/3}$.\\
First we want to show that the resolvent of $\CA_h\,$, as $h\rightarrow0$, can be conveniently approximated by the resolvent of operator
\beq
\left\{\begin{array}{l}
\CA_{-,h} = -h^2\frac{d^2}{dx^2}+i(V(0)+V'(0)x)\,,\\
\CD(\CA_{-,h}) = H_0^1(0,+\infty)\cap H^2(0,+\infty)\cap L^2(0,+\infty ; x^2dx)\,.
\end{array}\right.
\eeq
More precisely, given $\lambda_0>0$ and $\lambda(h) = \lambda_0 h^{2/3}$, we extend
$(\CA_h+\lambda_0h^{2/3})^{-1}$ on $(0,+\infty)$ by considering instead the operator
$\mathbf{1}_{[0,1]}(\CA_h+\lambda_0h^{2/3})^{-1}\mathbf{1}_{[0,1]}\,$. We then prove the following:
\begin{proposition}\label{propUpper1d}
Under the assumptions of Theorem \ref{thmNSAschro1_1d}, we have
 \beq\label{cvRes1}
 \|\mathbf{1}_{[0,1]}(\CA_h+\lambda_0h^{2/3})^{-1}\mathbf{1}_{[0,1]}-(\CA_{-,h}+\lambda_0h^{2/3})^{-1}\|_{\mathcal{L}(L^2(0,+\infty))} = 
 o\left(\frac{1}{h^{2/3}}\right)\,,
 \eeq
 as $h\rightarrow0\,$.
\end{proposition}
\textbf{Proof: }
Choosing again $\rho\in(1/3,2/3)\,$, we consider $\chi_{-,h} = \eta_{-,h}\mathbf{1}_{[0,1]}\,$, where
$\eta_{-,h}\in\mathcal{C}_0^\infty(-\infty,h^\rho ; [0,1])\,$, $\eta_{-,h}(x) = 1$ for $x\leq h^\rho/2\,$.
We set
\[
 \tilde\chi_h = \sqrt{1-\chi_{-,h}^2}\mathbf{1}_{[0,1]}\,,
\]
and we use the approximate resolvent
\[
 \tilde\CR(h) = \chi_{-,h}(\CA_{-,h}+\lambda_0h^{2/3})^{-1}\chi_{-,h}+\tilde\chi_h(\CA_h+\lambda_0h^{2/3})^{-1}\tilde\chi_h\,.
\]
We have
\begin{eqnarray*}
&& (\CA_h+\lambda_0h^{2/3})\tilde\CR(h) = 
I+\chi_{-,h}(\CA_h-\CA_{-,h})(\CA_{-,h}+\lambda_0h^{2/3})^{-1}\chi_{-,h} \\
&& +[\CA_h,\chi_{-,h}](\CA_{-,h}+\lambda_0h^{2/3})^{-1}\chi_{-,h}+
[\CA_h,\tilde\chi_h](\CA_h+\lambda_0h^{2/3})^{-1}\tilde\chi_h\,,
\end{eqnarray*}
hence, composing on the left by $\mathbf{1}_{[0,1]}(\CA_h+\lambda_0h^{2/3})^{-1}\mathbf{1}_{[0,1]}\,$,
\begin{eqnarray}
 && \mathbf{1}_{[0,1]}(\CA_h+\lambda_0h^{2/3})^{-1}\mathbf{1}_{[0,1]} - \chi_{-,h}(\CA_{-,h}+\lambda_0h^{2/3})^{-1}\chi_{-,h} =
 \tilde\chi_h(\CA_h+\lambda_0h^{2/3})^{-1}\tilde\chi_h \nonumber\\
 &&
 -\mathbf{1}_{[0,1]}(\CA_h+\lambda_0h^{2/3})^{-1}\chi_{-,h}(\CA_h-\CA_{-,h})(\CA_{-,h}+\lambda_0h^{2/3})^{-1}\chi_{-,h} \nonumber\\
 &&-\mathbf{1}_{[0,1]}(\CA_h+\lambda_0h^{2/3})^{-1}\mathbf{1}_{[0,1]}[\CA_h,\chi_{-,h}](\CA_{-,h}+\lambda_0h^{2/3})^{-1}\chi_{-,h} \nonumber\\
 &&-\mathbf{1}_{[0,1]}(\CA_h+\lambda_0h^{2/3})^{-1}\mathbf{1}_{[0,1]}[\CA_h,\tilde\chi_h](\CA_h+\lambda_0h^{2/3})^{-1}\tilde\chi_h\,.
 \label{decomTildeR}
\end{eqnarray}
We want to prove that the right-hand side behaves as $o(h^{-2/3})$ as $h\rightarrow 0\,$.\\
Consider first the second term. We have clearly
\beq\label{resAh+lambda}
 \|(\CA_h+\lambda_0h^{2/3})^{-1}\| = \mathcal{O}\left(\frac{1}{h^{2/3}}\right)\,,
\eeq
and
\beq\label{resA-h+lambda}
 \|(\CA_{-,h}+\lambda_0h^{2/3})^{-1}\| = \mathcal{O}\left(\frac{1}{h^{2/3}}\right)\,.
\eeq
Hence, we can easily check, as in (\ref{estRHS_R1}) after replacing $-\lambda(h)$ by $+\lambda_0h^{2/3}\,$, that
\beq\label{estResInter1b}
\|
\mathbf{1}_{[0,1]}(\CA_h+\lambda_0h^{2/3})^{-1}\chi_{-,h}(\CA_h-\CA_{-,h})(\CA_{-,h}+\lambda_0h^{2/3})^{-1}\chi_{-,h}
\| = \mathcal{O}\left(\frac{1}{h^{2(2/3-\rho)}}\right)\,.
\eeq
We can also check, as in (\ref{estCommut1}), that
\[
 \|[\CA_h,\chi_{-,h}](\CA_{-,h}+\lambda_0h^{2/3})^{-1}\chi_{-,h}\| = \mathcal{O}(h^{2/3-\rho})\,.
\]
Consequently, in view of (\ref{resAh+lambda}),
\beq\label{estCrochet1b}
\|\mathbf{1}_{[0,1]}(\CA_h+\lambda_0h^{2/3})^{-1}\mathbf{1}_{[0,1]}[\CA_h,\chi_{-,h}](\CA_{-,h}+\lambda_0h^{2/3})^{-1}\chi_{-,h}\| = 
\mathcal{O}\left(\frac{1}{h^\rho}\right)\,,
\eeq
and similarly,
\beq\label{estCrochet1bis}
\|\mathbf{1}_{[0,1]}(\CA_h+\lambda_0h^{2/3})^{-1}\mathbf{1}_{[0,1]}[\CA_h,\tilde\chi_h](\CA_h+\lambda_0h^{2/3})^{-1}\tilde\chi_h\| = 
\mathcal{O}\left(\frac{1}{h^\rho}\right)\,.
\eeq
Let us now consider the first term in the right-hand side of (\ref{decomTildeR}).
After replacing $\lambda_0h^{2/3}$ by $\lambda_0h^{2/3}-iV(0)$,
we can assume that $V(0) = 0$. By applying the equality
\[
|\Im\sca{(\CA_h+\lambda_0h^{2/3})v}{v}| = \|V^{1/2}v\|_{L^2(0,1)}^2\,,
\]
to $v = \tilde\chi_h(\CA_h+\lambda_0h^{2/3})^{-1}\tilde\chi_h u$, $u\in L^2(0,+\infty)$, and after noticing (since we assumed $V'>0$ on $[0,1]$)
that for some $C>0\,$,
\[
 \forall x\in\Supp\tilde\chi_h,~~~V(x)\geq V(h^\rho/2) \geq \frac{h^\rho}{C}\,,
\]
we get
\begin{eqnarray*}
&& \frac{h^\rho}{C}\|\tilde\chi_h(\CA_h+\lambda_0h^{2/3})^{-1}\tilde\chi_h u\|^2  \leq 
 \|V^{1/2}\tilde\chi_h(\CA_h+\lambda_0h^{2/3})^{-1}\tilde\chi_h u\|^2 \\
  &&\leq \|(\CA_h+\lambda_0h^{2/3})\tilde\chi_h(\CA_h+\lambda_0h^{2/3})^{-1}\tilde\chi_h u\|~
  \|\tilde\chi_h(\CA_h+\lambda_0h^{2/3})^{-1}\tilde\chi_h u\| \\
  & & \leq (\|u\| + \|[\CA_h,\tilde\chi_h](\CA_h+\lambda_0h^{2/3})^{-1}\tilde\chi_h u\|)\|\tilde\chi_h(\CA_h+\lambda_0h^{2/3})^{-1}\tilde\chi_h u\|\\
   & & \leq (1+\mathcal{O}(h^{2/3-\rho}))\|\tilde\chi_h(\CA_h+\lambda_0h^{2/3})^{-1}\tilde\chi_h u\|\|u\|\,,
\end{eqnarray*}
where we have used an estimate similar to (\ref{estCommut1}) (with $+\lambda_0h^{2/3}$ instead of $-\lambda(h)$) to control the
commutator term. \\
Hence,
\beq\label{estResInter1c}
\|\tilde\chi_h(\CA_h+\lambda_0h^{2/3})^{-1}\tilde\chi_h\| = \mathcal{O}\left(\frac{1}{h^\rho}\right)\,.
\eeq
Thus, since $\rho\in(1/3,2/3)$, (\ref{decomTildeR}), (\ref{estResInter1b}), (\ref{estCrochet1b}), (\ref{estCrochet1bis}) and
(\ref{estResInter1c}) yield
\beq\label{cvResInter1a}
\|\mathbf{1}_{[0,1]}(\CA_h+\lambda_0h^{2/3})^{-1}\mathbf{1}_{[0,1]} - \chi_{-,h}(\CA_{-,h}+\lambda_0h^{2/3})^{-1}\chi_{-,h}\| = 
o\left(\frac{1}{h^{2/3}}\right)\,.
\eeq
In order to get (\ref{cvRes1}), it remains to show that
\beq\label{cvResInter1b}
\|\chi_{-,h}(\CA_{-,h}+\lambda_0h^{2/3})^{-1}\chi_{-,h}-(\CA_{-,h}+\lambda_0h^{2/3})^{-1}\| = o\left(\frac{1}{h^{2/3}}\right)\,.
\eeq
In this purpose, we write
\[
 \chi_{-,h}(\CA_{-,h}+\lambda_0h^{2/3})^{-1}\chi_{-,h}(\CA_{-,h}+\lambda_0h^{2/3}) = \chi_{-,h}^2-
 \chi_{-,h}(\CA_{-,h}+\lambda_0h^{2/3})^{-1}[\CA_{-,h},\chi_{-,h}\,]\,,
\]
and composing on the right by $(\CA_{-,h}+\lambda_0h^{2/3})^{-1}$\,,
\begin{eqnarray*}
&& (\CA_{-,h}+\lambda_0h^{2/3})^{-1} - \chi_{-,h}(\CA_{-,h}+\lambda_0h^{2/3})^{-1}\chi_{-,h}  =  \tilde\chi_h^2(\CA_{-,h}+\lambda_0h^{2/3})^{-1}
\tilde\chi_h^2
\\ & &
+\chi_{-,h}(\CA_{-,h}+\lambda_0h^{2/3})^{-1}[\CA_{-,h},\chi_{-,h}](\CA_{-,h}+\lambda_0h^{2/3})^{-1}\,.
\end{eqnarray*}
The second term in the right-hand side can be estimated as (\ref{estCrochet1b}), while the first one satisfies the same bound as in
(\ref{estResInter1c}), hence (\ref{cvResInter1b}) holds. This concludes the proof of Proposition \ref{propUpper1d}.
\hfill $\square$\\

The upper bound (\ref{upper1d}) follows easily from \cite{Kato}, Section IV, \S$3.5$. Indeed, for any subsequence $h_j\rightarrow0$ and any
eigenvalue \[\mu\in h_j^{2/3}\sigma((\CA_{-,h_j}+\lambda_0h_j^{2/3})^{-1})\setminus\{0\}\,,\] there exists a sequence $(\mu_j)_{j\geq1}\,$, with
\[
 \mu_j\in h_j^{2/3}\sigma(\mathbf{1}_{[0,1]}(\CA_{h_j}+\lambda_0h_j^{2/3})^{-1}\mathbf{1}_{[0,1]})\setminus\{0\}
= h_j^{2/3}\sigma((\CA_{h_j}+\lambda_0h_j^{2/3})^{-1})\setminus\{0\}\,,
 \]
such that $\mu_j\rightarrow\mu$ as $j\rightarrow+\infty\,$. In particular, with $\mu = 1/(e^{i\pi/3}|\mu_1|J^{2/3}+\lambda_0)\,$, we get a sequence
$\lambda_j = h_j^{2/3}(1/\mu_j-\lambda_0)\in\sigma(\CA_{-,h_j})$ such that $h_j^{-2/3}\Re\lambda_j\rightarrow|\mu_1|J^{2/3}/2$ as
$j\rightarrow+\infty\,$, which proves (\ref{upper1d}).

\section{Upper bound for a Morse potential}\label{sthm3Upper}
In this section, we prove the upper bound in (\ref{limSpect2}), following the method of Section \ref{sthm2}.
Namely, we want to prove
\beq\label{upper2}
\varlimsup\limits_{h\to0}\frac{1}{h}\inf\Re\sigma(\CA_h) \leq \frac{\kappa}{2}\,.
\eeq
Let $\lambda_0>0$ and $\lambda(h) = \lambda_0h$. Let $k_0\in\{1,\dots,p\}$ such that 
\[\kappa_{k_0} = \min_{k=1,\dots,p}\kappa_k =: \kappa\,,\]
where $\kappa_k$ is the quantity defined in (\ref{defKappak}).\\
By reproducing the argument given at  the end of previous section, it is enough to prove the following:
\begin{proposition}\label{propUpper2}
Under the assumptions of Theorem \ref{thmNSAschro2}, we have
 \beq\label{cvRes2}
 \|\mathbf{1}_{\bar\Omega}(\CA_h+\lambda_0h)^{-1}\mathbf{1}_{\bar\Omega}-(\CH_{k_0,h}+\lambda_0h)^{-1}\|_{\mathcal{L}(L^2(0,+\infty))} = 
 o\left(\frac{1}{h}\right)\,,
 \eeq
as $h\rightarrow0$, where $\CH_{k_0,h}$ is the approximate operator defined in (\ref{defHkh}).
\end{proposition}
\textbf{Proof: }
Let us denote
\[
 \mathcal{L}(x_{k_0}^c) = V^{-1}(\{V(x_{k_0}^c)\})\setminus\{x_{k_0}^c\}\,.
\]
By assumption (\ref{hypNiveaux}), we have
\beq\label{noCritPt}
\forall x\in\mathcal{L}(x_{k_0}^c),~~~\nabla V(x) \neq 0\,.
\eeq
Notice that $V^{-1}(\{V(x_{k_0^c})\})$ may contain only $x_{k_0}^c$ if it is an absolute extremum.
Some points $x\in\mathcal{L}(x_{k_0}^c)$ could also lie on $\pa\Omega\,$, but to simplify the proof, we shall
assume that $\mathcal{L}(x_{k_0}^c)\cap\pa\Omega = \emptyset$. If not, we can handle the corresponding terms in (\ref{decomTildeQ}) by using estimates
similar to (\ref{estRHS_R1bdry}) and (\ref{estCommut1bdry}).\\
Let $\rho'\in(1/3,1/2)\,$, and $\theta_{k_0,h}\,$, $\CH_{k_0,h}$ as in Section \ref{sthm3Lower}. 
For every $h\in(0,h_0)\,$, we choose a set of indices $L(h)\subset\mathbb{N}$ and a set of points
\[
 \big\{d_\ell(h)\in\mathcal{L}(x_{k_0}^c) : \ell\in L(h)\big\}\,,
\]
such that
\[
 \mathcal{L}(x_{k_0}^c)\subset\bigcup_{\ell\in L(h)}B(d_\ell(h), h^{\rho'})
\]
and such that $\bar B(d_\ell(h),h^{\rho'}/2)\cap \bar B(d_m(h),h^{\rho'}/2) = \emptyset$ for every $\ell,m\in L(h)\,$, $\ell\neq m\,$.\\
For $\ell\in L(h)\,$, we choose 
$\varphi_{\ell,h}\in\mathcal{C}_0^\infty\big(B(d_\ell(h),h^{\rho'}) ; [0,1]\big)\,$, $\varphi_{\ell,h}(x) = 1$
if\\$x\in\bar B(d_\ell(h),h^{\rho'}/2)\,$. Let
\[
 \psi_h = \mathbf{1}_{\bar\Omega}(x)\sqrt{1-\theta_{k_0,h}(x)^2-\sum_{\ell\in L(h)}\varphi_{\ell,h}(x)^2}\,.
\]
We will use the same kind of approximate operator as before on $\Supp \varphi_{\ell,h}$:
\[
 \CA_{\ell,h} = -h^2\Delta+i\big(V(d_\ell(h))+\nabla V(d_\ell(h))\cdot(x-d_\ell(h))\big)\,.
\]
Let $\lambda_0>0\,$.
Our approximate resolvent is
\beq
\tilde{\mathcal{Q}}(h) = \theta_{k_0,h}(\CH_{k_0,h}+\lambda_0h)^{-1}\theta_{k_0,h} + \sum_{\ell\in L(h)}\varphi_{\ell,h}
(\CA_{\ell,h}+\lambda_0h)^{-1}\varphi_{\ell,h}
 + \psi_h(\CA_h+\lambda_0h)^{-1}\psi_h\,,\label{defQtilde}
\eeq
and we have
\begin{eqnarray}
 & & \mathbf{1}_{[0,1]}(\CA_h+\lambda_0h)^{-1}\mathbf{1}_{[0,1]} - \theta_{k_0,h}(\CH_{k_0,h}+\lambda_0h)^{-1}\theta_{k_0,h} =
 \psi_h(\CA_h+\lambda_0h)^{-1}\psi_h \nonumber\\
 & & + \sum_{\ell\in L(h)}\varphi_{\ell,h}(\CA_{\ell,h}+\lambda_0h)^{-1}\varphi_{\ell,h} \nonumber \\
 & & -\mathbf{1}_{\bar\Omega}(\CA_h+\lambda_0h)^{-1}\theta_{k_0,h}(\CA_h-\CH_{k_0,h})(\CH_{k_0,h}+\lambda_0h)^{-1}\theta_{k_0,h} \nonumber\\
 & &-\mathbf{1}_{\bar\Omega}(\CA_h+\lambda_0h)^{-1}\mathbf{1}_{\bar\Omega}[\CA_h,\theta_{k_0,h}](\CH_{k_0,h}+\lambda_0h)^{-1}\theta_{k_0,h} 
 \nonumber\\
 & & -\sum_{\ell\in L(h)}\mathbf{1}_{\bar\Omega}(\CA_h+\lambda_0h)^{-1}\varphi_{\ell,h}(\CA_h-\CA_{\ell,h})(\CA_{\ell,h}+\lambda_0h)^{-1}
 \varphi_{\ell,h} \nonumber \\
 & &-\sum_{\ell\in L(h)}\mathbf{1}_{\bar\Omega}(\CA_h+\lambda_0h)^{-1}\mathbf{1}_{\bar\Omega}[\CA_h,\varphi_{\ell,h}]
 (\CA_{\ell,h}+\lambda_0h)^{-1}\varphi_{\ell,h} 
 \nonumber \\
 & &-\mathbf{1}_{\bar\Omega}(\CA_h+\lambda_0h)^{-1}\mathbf{1}_{\bar\Omega}[\CA_h,\psi_h](\CA_h+\lambda_0h)^{-1}\psi_h\,.
 \label{decomTildeQ}
\end{eqnarray}
In the following, we assume for simplicity that $V(x_{k_0}^c) = 0$ (if not, one only has to replace $\lambda_0h$ by $\lambda_0h+iV(x_k^c)$)\,.\\
Using that, for some $C>0\,$,
\[
 \forall x\in\Supp\psi_h,~~~|V(x)|\geq\frac{h^{\rho'}}{C}\,,
\]
we can prove as (\ref{estResInter1c}) that
\[
 \|\psi_h(\CA_h+\lambda_0h)^{-1}\psi_h\| = \mathcal{O}\left(\frac{1}{h^{\rho'}}\right)\,.
\]
Besides, as already stated, we have by rescaling:
\beq\label{estresQuadUpper1}
\|(\CA_{\ell,h}+\lambda_0h)^{-1}\| = \mathcal{O}\left(\frac{1}{h^{2/3}}\right)~~~\textrm{ and }~~~
\|(\CH_{k_0,h}+\lambda_0h)^{-1}\| = \mathcal{O}\left(\frac{1}{h}\right)\,.
\eeq
We have also clearly
\beq\label{estResQuadUpper2}
\|(\CA_h+\lambda_0h)^{-1}\| = \mathcal{O}\left(\frac{1}{h}\right)\,.
\eeq
Hence, we can check, as in (\ref{estResInter1b}) and (\ref{estCrochet1b}), that all the other terms in the right-hand side of (\ref{decomTildeQ})
are of the form $o(h^{-1})$ as $h\rightarrow0$. The sums over $L(h)$ can be estimated by Lemma \ref{lemCotlar}.\\
Finally, it remains to show that 
\[
 \|\theta_{k_0,h}(\CH_{k_0,h}+\lambda_0h)^{-1}\theta_{k_0,h}-(\CH_{k_0,h}+\lambda_0h)^{-1}\| = o\left(\frac{1}{h}\right)\,,
\]
which can be done as for (\ref{cvResInter1b}).
\hfill $\square$

\section{Semigroups estimates}\label{sDecaySG}
In this section we prove Corollary \ref{corNSAschroSG} by using a quantitative version of the Gearhardt-Pr\"uss Theorem
\cite{HelSjo}. Indeed, the standard version does not enable us to get a uniform control of the constant $M_\varepsilon$ with
respect to $h$ in (\ref{estSG1}) and (\ref{estSG2}).\\
We focus on the proof of (\ref{estSG1}), the case of $(ii)$ being similar.\\
For all $\varepsilon>0\,$, according to (\ref{estRes1}) there exists $h_\varepsilon>0$ such that
\[
\forall h\in(0,h_\varepsilon),~~ \sup_{\tiny{\nu\in\mathbb{R}}}\|(\CA_h-(|\mu_1|J_m^{2/3}/2-\varepsilon)h^{2/3}-i\nu)^{-1}\|\leq \frac{C_\varepsilon}
{h^{2/3}}\,.
\]
Moreover, the operator $\CA_h$ being maximal accretive, it generates a contraction semigroup $e^{-t\CA_h}$: 
\beq\label{dissipAprioriSchro}
\forall t>0,~~~ \|e^{-t\CA_h}\|\leq 1\,.
\eeq
We apply \cite{HelSjo}, Theorem $1.5$, with
$\omega = -(|\mu_1|J_m^{2/3}-\varepsilon)h^{2/3}<0$, $r(\omega)^{-1}\leq C_\varepsilon h^{-2/3}\,$,
$m(t)\equiv1$ and $a=\tilde a = t/2$, which yields
\beq\label{dissipGPSchro}
\|e^{-t\CA_h}\|\leq \frac{(|\mu_1|J_m^{2/3}-\varepsilon)C_\varepsilon}{1-e^{-(|\mu_1|J_m^{2/3}/2-\varepsilon)h^{2/3}t/2}}
e^{-(|\mu_1|J_m^{2/3}/2-\varepsilon)h^{2/3}t}\,.
\eeq
Let $c_0>0$ and $t_h = 2c_0h^{-2/3}/(|\mu_1|J_m^{2/3}-\varepsilon)\,$.
Then, by (\ref{dissipGPSchro}),
\[
\forall t\geq t_h,~~~\|e^{-t\CA_h}\|\leq M_\varepsilon^{(1)}e^{-(|\mu_1|J_m^{2/3}-\varepsilon)h^{2/3}t}\,,
\]
with
\[
 M_\varepsilon^{(1)} = \frac{(|\mu_1|J_m^{2/3}-\varepsilon)C_\varepsilon}{1-e^{-c_0}}\,.
\]
Moreover, by (\ref{dissipAprioriSchro}),
\[
\forall t\leq t_h,~~~\|e^{-t\CA_h}\|\leq M_\varepsilon^{(2)}e^{-(|\mu_1|J_m^{2/3}-\varepsilon)h^{2/3}t}\,,
\]
with $M_\varepsilon^{(2)} = e^{2c_0}\,$. \\
Thus,
\beq\label{dissipinterSchro}
\forall t>0\,,~~~\|e^{-t\CA_h}\|\leq M_\varepsilon e^{-(|\mu_1|J_m^{2/3}-\varepsilon)h^{2/3}t}\,,
\eeq
with $M_\varepsilon = \max(M_\varepsilon^{(1)},M_\varepsilon^{(2)})$\,.\\
Estimate (\ref{estSG2}) can be proved the same way.\\

To prove the optimality statement in $(iii)$ of Corollary \ref{corNSAschroSG}, under the assumptions of Theorem \ref{thmNSAschro2},
we just consider
\[
 u_h\in\ker(\CA_h-\lambda_{0,h}h)\,,
\]
where $\lambda_{0,h}$ satisfies $h\lambda_{0,h}\in\sigma(\CA_h)$ and $h\Re\lambda_{0,h} = \inf\Re\sigma(\CA_h)$\,.\\
Then, we have
\[
 e^{-t\CA_h}u_h(x) = e^{-\lambda_{0,h}ht}u_h(x)\,.
\]
Thus, by (\ref{limSpect2}), for every $t>0$ and $\varepsilon>0\,$, there exists $h_\varepsilon>0$ such that, for every $h\in(0,h_\varepsilon)\,$,
\[
 \|e^{-t\CA_h}u_h\| = e^{-\lambda_{0,h}ht}\|u_h\|
  \geq e^{-(\kappa/2+\varepsilon)ht}\|u_h\|\,.
\]
Optimality in (\ref{estSG1}) under the assumptions of Theorem \ref{thmNSAschro1_1d} can be proved the same way.

\section{Application to the stability of the normal state in superconductivity}\label{sSuper}
\sectionmark{Application in superconductivity}
In this section, we recall the results of \cite{Alm} and explain how we can recover them, in the simplified setting of a smooth domain $\Omega\,$, by
rewriting Corollary \ref{corNSAschroSG} in the large domain limit.

\subsection{The time-dependent Ginzburg-Landau equations}
In this subsection we recall the time-dependent Ginzburg-Landau model, and we introduce the simplifications leading to the linear problem
which shall be considered in next subsection.\\

Superconducting materials are known to lose their electrical resistance when placed at a lower temperature than their critical one.
However, if a sufficiently strong current is applied throughout the sample, then superconductivity disappears and the material reverts to the
normal state, even if the temperature remains lower than the critical one.\\
In order to understand this phenomenon, we consider the time-dependent Ginzburg-Landau model, which can be written as
follows in the $2$-dimensional setting (see \cite{Alm} for the $3$-dimensional version of the system):
\begin{equation}\label{TDGL}
\left\{\begin{array}{ll}
        \pa_t\psi+i\Phi\psi=(\nabla-i\mathbf{A})^2\psi
        +\psi(1-|\psi|^2)\,, & (t,x)\in\mathbb{R}^+\times\Omega \\
\kappa^2\textrm{curl}^2\mathbf{A}+\sigma(\pa_t\mathbf{A}+\nabla\Phi)=\Im(\overline{\psi}(\nabla-i\mathbf{A})\psi)\,, &
(t,x)\in\mathbb{R}^+\times\Omega \\
\psi(t,x) = 0\,, & (t,x)\in\mathbb{R}^+\times\pa\Omega \\
\sigma(\pa_t\mathbf{A}(t,x)+\nabla\Phi(t,x))\cdot\vec n(x) = J(x), & (t,x)\in\mathbb{R}^+\times\pa\Omega \\
\textrm{curl}~\mathbf{A}(t,x) = H_{ex}(x)\,, & (t,x)\in\mathbb{R}^+\times\pa\Omega\,,\\
\psi(0,x) = \psi_0(x)\,, & x\in\Omega\,,\\
\mathbf{A}(0,x) = \mathbf{A}_0(x)\,,& x\in\Omega\,.
       \end{array}\right.
\end{equation}

Here $\Omega\subset\mathbb{R}^2$ is a smooth bounded, connected domain and $\vec n(x)$ denotes the 
outward normal on $\pa\Omega$ at $x$.\\
The unknown functions are $\psi(t,x)\in\mathbb{R}\,$, $\mathbf{A}(t,x)\in\mathbb{R}^2$ and $\Phi(t,x)\in\mathbb{R}\,$.
The function $\psi$ denotes the so-called
\emph{order parameter} of the superconductor, and $|\psi|^2$ represents the density of presence of superconducting electrons in the material.
Hence $\psi\equiv0$ corresponds to the normal state where superconductivity does not take place, whereas $\psi\equiv 1$ represents a purely
superconducting state.\\
$\mathbf{A}$ denotes the magnetic potential in the sample, and $\Phi$ the electric current. $H_{ex}$ denotes the exterior magnetic field, and
$J\in\mathcal{C}^2(\pa\Omega)$ represents the electric current applied through $\Omega\,$. In the following we will denote the magnetic field by
$\mathbf{B} = \textrm{curl }~\mathbf{A}\,$.\\
The constants $\kappa$ and $\sigma$ denote respectively the Ginzburg-Landau parameter, which is a material property, and the normal conductivity of
the sample.\\

Our goal in the following is to prove, under additional assumptions and in the large domain limit, that
if the applied electric current $J$ is strong enough, then the normal
state solution is stable as $t\rightarrow+\infty\,$. This problem was solved in \cite{Alm}, in a more physically relevant setting.
More precisely, the author considered a non-smooth domain $\Omega$ with right-angled corners at its boundary, such that
$\pa\Omega = \pa\Omega_c\sqcup\pa\Omega_i$, with different boundary conditions on each component $\pa\Omega_c$ and $\pa\Omega_i\,$. Here we shall
recover the results in \cite{Alm} in the case of a smooth boundary.\\

Let us first consider the stationary normal solution $(0,\mathbf{A}_n,\Phi_n)(x)$ of (\ref{TDGL}). Then the second line of (\ref{TDGL}) yields
\[
 \frac{\kappa^2}{\sigma}\textrm{curl}~\mathbf{B}_n + \nabla\Phi_n = 0\,,
\]
where $\mathbf{B}_n = \textrm{curl}~\mathbf{A}_n\,$. Hence $\Phi_n$ is harmonic in $\Omega\,$.\\
Now we neglect the effects of the magnetic field, that is, we assume $H_{ex} = \mathbf{A} = 0$, and we consider the linearization of
(\ref{TDGL}) near the stationary normal state $(0,0,\Phi_n)$ which leads to the system
\beq\label{linTDGL}
\left\{\begin{array}{ll}
       \pa_t\psi-\Delta\psi+i\Phi_n\psi-\psi = 0\,, & (t,x)\in\mathbb{R}^+\times\Omega\,, \\
       -\Delta\Phi_n = 0\,, & x\in\Omega, \\
       \nabla\phi_n(x)\cdot\vec n(x) = \frac{\kappa^2}{\sigma} J(x), & x\in\pa\Omega\,, \\
       \psi(t,x) = 0\,, & (t,x)\in\mathbb{R}^+\times\pa\Omega \\
       \psi(0,x) = \psi_0(x)\,, & x\in\Omega\,.\\
       \end{array}\right.
\eeq
Furthermore, we shall assume that $\nabla\Phi_n\neq0$ in $\Omega$. Indeed, in the setting of a domain $\Omega$ with right-angled corners at its
boundary, this assumption can be easily justified, see \cite{Alm}.\\

The case where the magnetic field is not neglected has been studied in \cite{AHP1,AHP2,AHP3,AH}. Moreover, the results of \cite{AH} include
the analysis of the nonlinear term $\psi(1-|\psi|^2)$.\\

In the following subsection, we shall assume that $\Omega$ is a large domain in order to recover the operator $\CA_h$ studied in the previous
sections.

\subsection{Stability of the normal state}
Here again we follow \cite{Alm}. We consider equations (\ref{linTDGL}) in the domain $\Omega_R = \{Rx : x\in\Omega\}$ for $R>1\,$. In order
to preserve the gradient, we consider an electric potentiel of the form
\[
 \Phi_R(x) = R\Phi_n\left(\frac{x}{R}\right)\,.
\]
Thus, we consider the problem
\beq\label{linTDGLR}
\left\{\begin{array}{ll}
       \pa_t\psi_R-\Delta\psi_R+iR\Phi_n\left(\frac{x}{R}\right)\psi_R-\psi_R = 0\,, & (t,x)\in\mathbb{R}^+\times\Omega_R\,, \\
       \psi_R(t,x) = 0\,, & (t,x)\in\mathbb{R}^+\times\pa\Omega \\
       \psi_R(0,x) = \psi_{0,R}(x)\,, & x\in\Omega\,,\\
       \end{array}\right.
\eeq
where $\Omega\subset\mathbb{R}^n$ and we no longer need to assume $n=2$.\\
The function $\Phi_n$ is assumed to be smooth and satisfies, for all $x\in\bar\Omega\,$, $\nabla\Phi_n(x)\neq0$\,.\\
In other words, we have
\beq\label{evolLR}
 \forall (t,x)\in\mathbb{R}^+\times\Omega,~~~\psi_R(t,x) = e^{-t\CL_R}\psi_{0,R}(x)\,,
\eeq
where
\[
\left\{ \begin{array}{l}
         \CL_R = -\Delta+iR\Phi_n\left(\frac{x}{R}\right)-1\,,\\
\CD(\CL_R) = H_0^1(\Omega)\cap H^2(\Omega)\,.
         \end{array}\right.
\]
Let us set $h = h(R) = R^{-3/2}\,,$ and $T_R : u(x) \mapsto Ru(x/R)\,$. Then,
\[
 T_R\CL_RT_R^{-1} = h(R)^{-2/3}(-h(R)^2\Delta +i\Phi_n-h(R)^{2/3}))\,.
\]
Hence, for all $t>0$ and $R>1$ we have $\|e^{-t\CL_R}\| = \|e^{-(th(R)^{-2/3})(\CA_{h(R)}-h(R)^{2/3})}\|$\,, where $\CA_h$ is the operator
defined in (\ref{defAh}), where $V = \Phi_n$ satisfies the assumptions of Theorem \ref{thmNSAschro1}. Thus, Corollary \ref{corNSAschroSG} yields:

\begin{theorem}[Y. Almog, \cite{Alm}]
Let
\[
 \pa\Omega_\perp = \{x\in\pa\Omega : \nabla\Phi_n(x)\times\vec n(x) = 0\}
\]
and, if $\pa\Omega_\perp\neq\emptyset\,$,
\[
 J_m = \min_{x\in\pa\Omega_\perp}|\nabla\Phi_n(x)|\,.
\]
Let
\[
 J_c = \left(\frac{2}{|\mu_1|}\right)^{3/2}\,.
\]
If $\pa\Omega_\perp = \emptyset$ or $J_m>J_c\,$, then for all $\varepsilon>0$, there exists $R_0>1$ and $M_\varepsilon>0$ 
such that, for all $R\geq R_0$ and $\psi_{0,R}\in H_0^1(\Omega)\cap
H^2(\Omega)\,$, the solution $\psi_R$ of (\ref{linTDGLR}) satisfies
\beq
\forall t>0\,,~~~\|\psi_R(t,\cdot)\|_{L^2(\Omega)} \leq M_\varepsilon\exp(-((J_m/J_c)^{2/3}-1-\varepsilon)t)\|\psi_{0,R}\|_{L^2(\Omega)}\,.
\eeq
\end{theorem}
The results of \cite{AHP1,AHP2,AHP3,AH} give similar conditions for the stability of the normal state in the presence of a magnetic field.

\end{document}